\documentclass[10pt]{article}

\usepackage{amsmath,amssymb,amsthm}
\usepackage{algorithm,algorithmic}
\usepackage{graphicx}
\usepackage{hyperref,cleveref}
\input{macros.sty}
\usepackage[textsize=small]{todonotes}

\usepackage[margin=1.1in]{geometry}
\newtheorem{definition}{Definition}

\usepackage{caption,color}

\title{Efficient algorithms for Bayesian Inverse Problems with Whittle--Mat\'ern Priors\thanks{{HA is partially supported by National Science Foundation (NSF) grants DMS-2110263, DMS-1913004 and the Air Force Office of Scientific Research under Award NO: FA9550-22-1-0248. AKS is partially supported by NSF DMS-1745654 and DMS-2026830.}}}
\author{Harbir Antil\thanks{{Department of Mathematical Sciences and the Center for Mathematics and Artificial Intelligence (CMAI),} George Mason University, USA. {hantil@gmu.edu}. ORCID: 0000-0002-6641-1449}  \and Arvind K.\ Saibaba\thanks{Department of Mathematics, North Carolina State University, USA. {asaibab@ncsu.edu}. ORCID: 0000-0002-8698-6100}}

\begin{document}

\maketitle

\begin{abstract}
This paper tackles efficient methods for Bayesian inverse problems with priors based on Whittle--Mat\'ern Gaussian random fields. The Whittle--Mat\'ern prior is characterized by a mean function and a covariance operator that is taken as a negative power of an elliptic differential operator. This approach is flexible in that it can incorporate a wide range of prior information including non-stationary effects, but it is currently computationally advantageous only for integer values of the exponent. In this paper, we derive an efficient method for handling all admissible noninteger values of the exponent. The method first discretizes the covariance operator using finite elements and quadrature, and uses preconditioned Krylov subspace solvers for shifted linear systems to efficiently apply the resulting covariance matrix to a vector. This approach can be used for generating samples from the distribution in two different ways: by solving a stochastic partial differential equation, and by using a truncated Karhunen-Lo\`eve expansion. We show how to incorporate this prior representation into the infinite-dimensional Bayesian formulation, and show how to efficiently compute the maximum a posteriori estimate, and approximate the posterior variance. Although the focus of this paper is on Bayesian inverse problems, the techniques developed here are applicable to solving systems with fractional Laplacians and Gaussian random fields. Numerical experiments demonstrate the performance and scalability of the solvers and their applicability to model and real-data inverse problems in tomography and a time-dependent heat equation. 
\end{abstract}


\section{Motivation and Introduction}
Inverse problems involve the use of experimental data or measurements to make inferences about parameters (e.g., initial or boundary conditions), governing the model. Inverse problems have a wide range of applications such as in material science, geophysics, and medicine. The Bayesian approach to inverse problems is prevalent because of its ability to incorporate prior knowledge of the parameters and its ability to quantify the uncertainty associated with the parameter reconstructions. However, many computational challenges persist in the implementation of the Bayesian approach to large-scale inverse problems.

Stuart~\cite{stuart2010inverse} advocated the use of Gaussian priors within Bayesian inverse problems with covariance operators of the form $\mc{C} = \mc{A}^{-\alpha}$ where $\mc{A}$ is an elliptic differential operator, the exponent $\alpha > d/2 $, and $d$ is the spatial dimension. This choice ensures that the covariance operator is trace-class and under appropriate conditions on the likelihood, the resulting inverse problem is well-posed. A computational framework for infinite-dimensional Bayes inverse problems based on Stuart's framework was developed in~\cite{bui2013computational}. For computational convenience, the authors used the covariance operator $\mc{C} = \mc{A}^{-2}$ since it satisfies the requirement $\alpha > d/2$ for up to three spatial dimensions. This choice has two advantages: first, it is integer-valued, and second, it is easy to obtain a factored form of the covariance operator. For these reasons, this approach is attractive from a computational standpoint, and has since become popular and has resulted in scalable software implementations~\cite{villa2018hippylib,kim2021hippylib}. A similar approach can also be found in~\cite{roininen2014whittle}, in which the authors also used the Whittle--Mat\'ern priors for Bayesian inverse problems. However, the choice of $\alpha = 2$ for computational convenience makes it restrictive in practice, and in many instances, it can lead to over-smoothing of reconstructions. It is, therefore, desirable to have an efficient computational method for handling exponents $\alpha > d/2$ that allows the user to choose the appropriate prior covariance based on prior knowledge or estimate it from data in a hierarchical Bayesian setting.

The use of Whittle--Mat\'ern priors in Bayesian inverse problems has close connections with the developments in Gaussian random fields in spatial and computational statistics. In~\cite{lindgren2011explicit}, the authors considered the stochastic partial differential equation (SPDE) approach to Gaussian random fields. They showed an explicit connection between Gaussian and Gauss-Markov random fields based on the Mat\'ern covariance family and provide explicit expressions for the precision operator for integer values of $\alpha$. Based on this connection, they provided extensions beyond the Mat\'ern model and showed how the SPDE approach can be applied on manifolds, for nonstationary random fields. In a review paper~\cite{lindgren2022spde}, 10 years since the publication of~\cite{lindgren2011explicit}, the authors trace the recent developments in the SPDE approach to Gaussian and non-Gaussian random fields. Recent work (e.g.,~\cite{bolin2018weak,bolin2020numerical}) allows for noninteger values of the exponent; see discussion below.

\paragraph{Goals and Main Contributions} In this paper, we want to develop an efficient method for representing and computing with Whittle--Mat\'ern Gaussian random priors in which the covariance operator is of the form $\mc{C} = \mc{A}^{-\alpha}$ and $\alpha > d/2$. Furthermore, we also want to incorporate this efficient representation of covariance operators while solving Bayesian inverse problems. The novel and noteworthy features of our contributions are:
\begin{enumerate}
    \item  We consider a quadrature based approach for the approximation of the covariance operator for fractional powers. The resulting approximation has the form of a rational function. The action of the covariance operator on a vector can be written as a sequence of shifted linear systems, which we accelerate using Krylov subspace methods. As a result, we have a method for efficiently applying the covariance operator for any exponent $\alpha > d/2$. The method is scalable to large problems and the cost is nearly the same for values of $\alpha$ between consecutive integers. This approach is applicable for stationary (including isotropic and anisotropic) and non-stationary random fields. 
    \item We leverage the fast method for covariance matrices to efficiently generate samples in two different ways: by solving the underlying SPDE that defines the prior covariance and by efficiently {computing} the truncated Karhunen-Lo\'eve approximation. 

    \item We show how to use this efficient representation of the prior covariance operator in an infinite-dimensional Bayesian problem formulation. We perform careful discretization of all the relevant quantities and leverage the fast applications of the covariance operator to compute the maximum a posteriori (MAP) estimate  and approximate posterior variance. 

    \item Although the focus of this paper is on Bayesian inverse problems, the computational techniques developed here are applicable to solving systems with fractional Laplacians \cite{bonito2015numerical,CJWeiss_BGVBWaanders_HAntil_2020a} and Gaussian random fields. 
\end{enumerate}
Through numerical experiments, we demonstrate the computational performance of the solvers including a $\sim 40$X speedup over na\"ive approaches.  We also demonstrate the feasibility of the approach on model and real-data problems from seismic and X-ray tomography, and an application to a 2D PDE-based inverse diffusion equation. 

\paragraph{Related Work} As mentioned earlier, in Bayesian inverse problems, the exponent is typically chosen to be $\alpha = 2$. To our knowledge, in Bayesian inverse problems, the problem of efficient computations with covariance operators for non-integer values of $\alpha$ is still an {outstanding} challenge. The fractional Laplacian with $0 < \alpha < 1$ has been used as a regularization operator in a non-Bayesian {setting~\cite{HAntil_SBartels_2017a,antil2020bilevel}. Notice that, \cite{HAntil_SBartels_2017a} uses the Fourier approach which is limited to periodic settings and \cite{antil2020bilevel}} uses an eigendecomposition of the Laplacian that is not scalable to large problems.

There are some techniques in computational and spatial statistics that are relevant to this discussion{; however for a detailed review, see~\cite{lindgren2022spde}}. In~\cite{simpson2008krylov}, a method was proposed for sampling from generalized Mat\'ern fields on compact Riemannian manifolds using the so-called matrix transfer technique. In this method, after discretizing the differential operator, the samples are obtained by applying the matrix function ($f(x) = x^{-\alpha}$) on a random vector. The application of the matrix function is accomplished using contour integral and Krylov subspace methods similar to ours. But our approach differs in how we approximate the fractional operator and also in the choice of Krylov subspace methods. Recent work in~\cite{bolin2018weak,bolin2020numerical} develop efficient methods for sampling from the SPDE for certain fractional exponents (in our notation, $d/4 < \alpha/2 < 1$); these references also use the same quadrature scheme for the integral representation of the fractional inverse as ours. However, our approach generalizes it to all exponents $\alpha > d/2$ and the use of Krylov subspace solvers substantially speeds up the solution of the SPDE. In a follow-up work~\cite{bolin2020rational}, the authors developed the rational SPDE approach in which they use a rational approximation of the form (in our notation) $\mc{A}^{-\alpha} \approx p(\mc{A})q(\mc{A})^{-1}$, where $p,q$ are polynomials. When the degrees of the polynomials are small, the discretized representations retain the computational benefits of the integer versions. However, the error analysis in~\cite[Remark 3.4]{bolin2020rational} suggests in order to balance the errors due to the spatial discretization and the error in the rational approximation, the degrees of the polynomials $p,q$ may be large, which leads to a loss in computational efficiency. {The paper in~\cite{herrmann2020multilevel} also works in a similar setting as ours in bounded domains but solves the linear systems in parallel using multilevel preconditioned iterative solvers with linear complexity in the degrees of freedom; additionally, their approach is also applicable to compact metric spaces.} Also relevant to this discussion are~\cite{lang2021galerkin,jansson2022surface} which deal with sampling from random fields on compact Riemannian manifolds. Our efficient approach for representing covariance matrices and sampling from the random field may be applicable to that setting as well. To our knowledge, these techniques have not been used within the context of Bayesian inverse problems, which is the primary focus of this paper.

\section{Background} In Section~\ref{ssec:prior}, we review the SPDE approach to Gaussian random fields  which form the basis for the prior distributions in Bayesian inverse problems. In Section~\ref{ssec:fractional}, we review the appropriate definitions of fractional elliptic operators and show how to discretize them in Section~\ref{ssec:disc}.

\subsection{SPDE approach to Whittle--Mat\'ern Gaussian Random Fields}\label{ssec:prior} 
We follow the Whittle--Mat\'ern {random field} approach, which models the random fields as stationary Gaussian random fields with the covariance kernel
\[ c(\bfx,\bfy) := \sigma^2 \frac{2^{1-\nu}}{\Gamma(\nu)}\left( \kappa {\|\bfx-\bfy\|_2 }\right)^\nu K_\nu\left(\kappa\|\bfx-\bfy\|_2 \right) , \qquad \bfx, \bfy \in \R^d,\]
where $\nu > 0$ is the smoothness parameter, the integer value of which controls the mean-square differentiability of the process, $K_\nu$ is the modified Bessel function of the second kind of order $\nu$, and $\kappa$ is a scaling parameter that governs the length scale. When $\nu = 1/2$, this corresponds to the exponential covariance function, and as $\nu\rightarrow \infty$, this corresponds to the Gaussian kernel. On $\R^d$, samples from the random field with the covariance function correspond to the solution of the stochastic partial differential equation (SPDE)
\begin{equation}\label{eqn:spde} (\kappa^2  - \Delta)^{(\nu + d/2)/2} u(\bfx) = W(\bfx),  \qquad \bfx \in \R^d, \end{equation}
where $W$ is a spatial Gaussian white noise process with unit variance on $\R^d$,  and the marginal variance is 
\[ \sigma^2 = \frac{\Gamma(\nu)}{(4\pi)^{d/2}\Gamma(\nu + d/2)}.\]
However, to implement this on a bounded region $\Omega \subset \R^d$, we need to impose additional boundary conditions. The choice of the boundary conditions may affect the correlation near the boundaries. To mitigate this effect, one may solve the SPDE on a slightly larger domain that encloses $\Omega$; other approaches have been proposed in~\cite{roininen2014whittle,daon2018mitigating}. The covariance operator corresponding to the  covariance kernel $c(\bfx,\bfy)$ is $ (\kappa^2  - \Delta)^{-(\nu + d/2)}$. In this paper, work with covariance operators of the form $\mc{C}_{\alpha}$
\begin{equation}\label{eqn:coviso} \mc{C}_{\alpha} := (\kappa^2  - \Delta)^{-\alpha}, \end{equation}
and  $\alpha :=\nu+ d/2$ with corresponding zero Neumann boundary conditions. The definition  ensures that the covariance operator $\mc{C}$ is a self-adjoint, trace-class operator~\cite{stuart2010inverse}.

Furthermore, since it more convenient to work with integer values of the exponent $\alpha$, the constraint $\alpha > d/2$ forces us to choose $\alpha = 1$ for one spatial dimension and $\alpha = 2 $ for two or three spatial dimensions. However, in the context of Bayesian inverse problems this choice of $\alpha$ can result in oversmoothing of the solution and  poor edge-preserving behavior of the solution. Furthermore, it is common to use parameters such as $\nu = 1/2,3/2,5/2$, etc; for example $\nu = 1/2$ corresponds to the well-known exponential covariance kernel. However, when $d=2$, $\alpha = 3/2$ which cannot be tackled using standard approaches. It is desirable to develop a numerical approach that can exploit the full power of these parameterized priors, without the constraints of integer powers. This motivates the use of fractional operators which we discuss in Section~\ref{ssec:fractional}.

Following~\cite[Section 3.2]{lindgren2011explicit} and~\cite[Section 2.6]{lindgren2022spde}, It is also straightforward to model non-stationary effects by considering the SPDE
\begin{equation}\label{eqn:spde1} 
(\kappa^2(\bfx)  - \nabla\cdot (\bfTheta (\bfx) \nabla ) )^{\alpha/2}u = \tau(\bfx) W, \qquad \bfx \in \R^d,
\end{equation}
where $\kappa^2(\bfx)$, $\tau(\bfx)$, and $\bfTheta(\bfx)$ are now spatially varying functions and $W$ is a spatial Gaussian white noise process with unit variance. Assuming that $\tau = 1$, the corresponding covariance operator (with zero Neumann boundary conditions) is 
\begin{equation}\label{eqn:covalpha} C_{\alpha} := ( \kappa^2(\bfx)  - \nabla\cdot (\bfTheta (\bfx) \nabla ))^{-\alpha}.  \end{equation}
Note that this assumes~\eqref{eqn:coviso} as a special case.

\subsection{Fractional Laplacian}\label{ssec:fractional}

{
For $\kappa^2 > 0$, let us use the notation $\kappa^2 - \Delta$ to denote the realization in $L^2(\Omega)$ of $\kappa^2 - \Delta$ with zero Neumann boundary conditions. Then there exists a sequence of eigenvalues $\{\lambda_k\}_{k\ge 1}$ satisfying
$0 < \lambda_1 \le \lambda_2 \le \cdots \le \lambda_k \le \cdots$ with $\lim_{k\rightarrow \infty} \lambda_k = \infty$.  The corresponding eigenfunctions $\{\varphi_k\}_{k\ge 1} \subset H^1(\Omega)$. It is well-known that, $\{\varphi_k\}_{k\ge 1}$ forms the orthonormal basis of $L^2(\Omega)$. Towards, this end, for any $s \ge 0$, we can define the fractional order Sobolev space $H^s(\Omega)$ as
\[
 H^s(\Omega) 
        := \left\{ u = \sum_{k=1}^\infty u_k \varphi_k \ : \  \| u \|_{H^s(\Omega)}^2
            := \sum_{k=1}^\infty \lambda_k^2 u_k^2 < \infty \right\} .
\]
We are now ready to define the fractional powers of $\kappa^2 - \Delta$, see \cite[Section~8]{HAntil_JPfefferer_SRogovs_2018a}. 
\begin{definition}
    Let $u \in C^\infty(\overline\Omega)$ with $\bfn \cdot \nabla u = 0$ on $\partial\Omega$, where $\bfn$ is the outward normal vector. 
    Then the fractional power of $\kappa^2 - \Delta$ is given by
    \[
        (\kappa^2 - \Delta)^s u := \sum_{k=1}^\infty \lambda_k^s u_k \varphi_k 
        \quad \mbox{with} \quad 
        u_k = \int_\Omega u \varphi_k . 
    \]
\end{definition}
Notice that $(\kappa^2 - \Delta)^s$ can be extended to an operator mapping from $H^s(\Omega)$ to $H^{-s}(\Omega)$ where $H^{-s}(\Omega)$ is the dual of $H^s(\Omega)$. We further emphasize that the condition $\kappa^2 > 0$ can be dropped in the above definition, see \cite[Def.~2.2]{HAntil_JPfefferer_SRogovs_2018a}. We also refer to \cite{HAntil_JPfefferer_SRogovs_2018a}, for extensions to nonzero Neumann conditions.

Finally, we can extend the above definition to $( \kappa^2(\bfx)  - \nabla\cdot (\bfTheta (\bfx) \nabla ))$ in \eqref{eqn:spde1}. In this case, we need to assume $\kappa^2 \in L^\infty(\Omega)$ and $\bfTheta$ to be symmetric, bounded in $L^\infty(\Omega)$ and elliptic.
}

\subsection{Discretization} \label{ssec:disc}

We show how to discretize the appropriate quantities in a finite dimensional space. The discussion here closely mimics the derivation in~\cite{bui2013computational}.
\paragraph{Finite dimensional parameter space}  Let $\phi_1,\dots,\phi_n$ be a basis for the finite dimensional space $\mc{V}_h \subset L^2(\Omega)$, with $N_h = \text{dim}(\mc{V}_h)$ where the subscript $h$ refers to a mesh discretization parameter. Assume that the basis functions corresponding to the nodal points $\{\bfx_j\}_{j=1}^{N_h}$ satisfy $\phi_j(\bfx_i) =\delta_{ij}$ for $1 \leq i, j \leq N_h$. We can then approximate the inversion parameter within this finite dimensional subspace as $m_h = \sum_{j=1}^{N_h}m_{h,j}\phi_j$, and we denote the vector of coefficients in the expansion $\bfm = \bmat{m_{h,1} & \dots & m_{h,N_h} }^\top \in \R^{N_h}.$

For two functions $u,v \in L^2(\Omega)$, the $L^2$ inner product  can be approximated as  $(u,v)_{L^2(\Omega)} \approx \langle \bfu, \bfv\rangle_{\bfM} = \bfu^\top \bfM \bfv$, where  $\bfM$ is the mass matrix with entries 
\[ \bfM_{ij} = \int_\Omega \phi_i(\bfx) \phi_j(\bfx) d\bfx \qquad 1 \leq i, j\leq N_h.\]
We denote $\R^{N_h}_{\bfM}$ to be the vector space $\R^{N_h}$ with the inner product $\langle \cdot, \cdot\rangle_{\bfM}$, to distinguish from the usual Euclidean space.

Let $\mc{B} : L^2(\Omega) \rightarrow L^2(\Omega)$. Then the matrix approximation of $\mc{B}_h : \mc{V}_h \rightarrow \mc{V}_h$, denoted $\bfB : \R^{N_h}_{\bfM} \rightarrow \R^{N_h}_{\bfM}$, is obtained as follows. First, we construct a matrix $\bf{S}$ with entries
\[ \bfS_{ij} = \int_\Omega  \phi_i(\bfx)\mc{B} \phi_j(\bfx) d\bfx = \langle\bfe_i, \bfB \bfe_j \rangle_{\bfM}\qquad 1 \leq i, j\leq N_h.\]
Here $\phi_i$ are the finite dimensional basis vectors and $\bfe_i$ is the canonical basis vector for $\R^{N_h}$ corresponding to the basis function $\phi_i$. This gives us the matrix representation $\bfB = \bfM^{-1}\bfS$. The matrix transpose of $\bfB$, denoted $\bfB^T$, has entries $[\bfB^T]_{ij} = b_{ji}$, but the adjoint of $\bfB : \R^{N_h}_{\bfM} \rightarrow \R^{N_h}_{\bfM}$, denoted $\bfB^*$, satisfies $\bfB^* = \bfM^{-1}\bfB^T \bfM$. 

\paragraph{Discretization of the Covariance operator}
In the following discussion, we first assume that $0 < s < 1$ and $\mc{A} :=  \kappa^2 - \nabla\cdot (\bfTheta(\bfx) \nabla)$, with zero Neumann boundary conditions.  To apply the fractional operator $u = \mc{A}^{-s}f$, we use sinc quadrature to approximate the integral~\cite{bonito2015numerical}. Let $z_{h,j} \in \mc{V}_h $ solve 
\[ \int_\Omega \left( (\bfTheta \nabla z_{h,j}) \cdot \nabla v_h +   (\kappa^2 + e^{j\zeta }) z_{h,j} v_h \right) dx  = \int_\Omega fv_h dx,  \quad \forall  v_h \in \mc{V}_h, -M_- \leq j \leq M_+,\]
and compute the approximation to $u$ as 
\[ u_h = \frac{\zeta\sin(s \pi)}{\pi} \sum_{j=-{M_{-}}}^{M_+}e^{(1-s)j\zeta } z_{h,j}, \]
where $M_+ = \lceil \frac{\pi^2}{4s\zeta^2} \rceil, 
M_- = \lceil \frac{\pi^2}{4(1-s)\zeta^2} \rceil$, and $\zeta = 1/\log(1/h)$. The number of terms in the expansion was chosen to balance the errors in the quadrature with the finite element discretization error. Alternatively, in matrix notation we can write 
\[\bfu= \sum_{j=-M_-}^{M_+} w_j ( \bfK + z_j\bfM)^{-1} \bfM\bff, \]
where $\bfK$ is the stiffness matrix with entries $$[\bfK]_{ij} = \int_\Omega  (\bfTheta \nabla \phi_i) \cdot \nabla \phi_j + \kappa^2 \phi_i \phi_j d\bfx \qquad 1 \leq i,j \leq N_h$$ and $\bfM$ is the mass matrix with entries $m_{ij} = \int_{\Omega} \phi_i\phi_jd\bfx$ for $1\leq i,j\leq N_h$. The quadrature weights and nodes are  
\begin{equation}\label{eqn:inverseweights}w_j = \frac{\zeta\sin(\alpha \pi)}{\pi} e^{(1-s)j\zeta }, \qquad  z_j = e^{j\zeta}, \quad -M_- \leq j \leq M_+. \end{equation}
Based on the preceding discussion, the matrix representation of the discretized version of $\mc{D} = \mc{A}^{-s}$ is 
\begin{equation*}\bfD := \sum_{j=-M_-}^{M_+} w_j   ( \bfK + z_j\bfM)^{-1}\bfM.
\end{equation*}
Let $\mathbb{N}$ denote the set of natural numbers (excluding zero). Now consider the covariance operator $\mc{C}_{\alpha} = \mc{A}^{-\alpha}$ where $\alpha > d/2$ and $\alpha \notin \mathbb{N}$.  We write $\alpha = r + s$ where $r = \lfloor\alpha \rfloor$ be the integer part of $\alpha$ and $0 < s < 1$ is the fractional part. Following~\cite[Remark 3.6]{bolin2020numerical}, the discretized representation of the covariance operator $\mc{C}_{\alpha} = \mc{A}^{-(\nu+d/2)}$ is
\begin{equation}\label{eqn:discov}\bfC_{\alpha} := \left(\sum_{j=-M_-}^{M_+} w_j   ( \bfK + z_j\bfM)^{-1}\bfM\right) \left[\bfK^{-1} \bfM\right]^{r}.
\end{equation}
Note that we can interchange the order of the fractional and the integer part since terms like $( \bfK + z_j\bfM)^{-1}\bfM$ and $\bfK^{-1} \bfM$ commute. Furthermore, the matrix $\bfC_\alpha$ is self-adjoint, i.e., $\bfC_\alpha^* =\bfC_\alpha$. An important point to note here is that $\bfC_{\alpha}$ is not formed explicitly. The application of the matrix to a vector can be obtained by $M_- + M_++1$ independent linear system solves which can be readily parallelized. However, we will show {how to compute} the action of $\bfC_{\alpha}$ efficiently using Krylov subspace methods for shifted linear systems (see Section~\ref{ssec:applying}).

\section{Efficient computations with the covariance operator}
In this section, we review the Multipreconditioned GMRES method for shifted linear systems~\cite{bakhos2017multipreconditioned} (MPGMRES-Sh) for solving sequences for shifted linear systems (Section~\ref{ssec:mpgmres}). We use this solver to accelerate the computations with the covariance operator $\bfC_{\alpha}$ and generating samples from the SPDE (Section~\ref{ssec:applying}). Finally, in Section~\ref{ssec:kle}, we also describe how to efficiently compute the truncated Karhunen-Lo\`eve expansion.
\subsection{Krylov methods for shifted systems}\label{ssec:mpgmres} The dominant cost of applying the fractional operator is the solution of a system of shifted linear equations. 
We first briefly review the MPGMRES-Sh approach~\cite{bakhos2017multipreconditioned}. This is an efficient method for solving shifted linear systems of the form 
\[ (\bfA_1 + \sigma_j\bfA_2)\bfx_j = \bfd \qquad 1 \leq  j \leq N_\sigma, \]
where $\bfA_1, \bfA_2 \in \R^{n\times n}$ and $\{\sigma_j\}_{j=1}^{N_\sigma}$ are a set of shifts. This solver can be used in multiple ways to accelerate computations involving Gaussian random fields.

In the MPGMRES-Sh approach, we select a set of shifts $\{\tau_j\}_{j=1}^{n_p}$ that determine the preconditioners $\bfP_j = (\bfA_1 + \tau_j \bfA_2).$ The method uses multiple preconditioners to build a basis $\bfZ_m \in \R^{n \times mn_p}$ that we will use to define a search space for all the shifted systems. In the below discussion, for simplicity, we drop the index on $\sigma$. At the $k$-th iteration, given the matrix $\bfV^{(k)} \in \R^{n\times n_p^{k-1}}$ the method builds the matrix $\bfZ^{(k)}$ as 
\[ \bfZ^{(k)} =  \bmat{\bfP_1^{-1}\widehat{\bfv}_k & \dots & \bfP_{n_p}^{-1} \widehat{\bfv}_k} \in \R^{n\times n_p} \qquad \bfv_k = \bfV^{(k)}\bfe_{n_p}.\]
The method is initialized with $\bfV^{(1)} = \bfd/\|\bfd\|_2$. We define the matrices that collect the columns $\bfV^{(k)}$ and $\bfZ^{(k)}$ as $\bfV_m = \bmat{\bfV^{(1)} & \dots & \bfV^{(m)}}$ and  $\bfZ_m = \bmat{\bfZ^{(1)} & \dots & \bfZ^{(m)}}$.  
Using these relationship we can derive the multipreconditioned Arnoldi relationship for a shift $\sigma$
\[ (\bfA_1 + \sigma \bfA_2 )\bfZ_m = \bfV_m \left( \bmat{\bfE_m \\ \bfzero} + \bar{\bfH}_m (\sigma \bfI - \bfT_m)\right) = \bfV_{m+1}\bar{\bfH}_m(\sigma;\bfT_m).\]
Here the matrices $\bfT_m$ and $\bfE_m$ are defined as follows. Let $\otimes$ denote the the Kronecker product and define
\[ \bfT^{(k)} = \blkdiag(\bfI_{n_p^{k-1} } \otimes \tau_1, \dots , \bfI_{n_p^{k-1} } \otimes \tau_{n_p}) \in \R^{n_p^k\times n_p^k}\]
and 
$\bfE^{(k)} = \bfe_{n_p}^\top \otimes \bfI_{n_p^{k-1}}$ for $1 \leq k \leq m$. Here $\blkdiag$ means a block diagonal matrix composed of the given subblocks. With these definitions in place, then $\bfT_m = \blkdiag(\bfT^{(1)}, \dots, \bfT^{(m)})$ and $\bfE_m = \blkdiag(\bfE^{(1)},\dots, \bfE^{(m)}).$ The matrix $\bfH_m$ is block upper-Hessenberg.

The main point here is that the solution for each shifted system corresponding to shift $\sigma$ is obtained as $\bfx_m(\sigma) = \bfZ_m\bfy_m(\sigma)$, where $\bfy_m (\sigma)\in \R^{mn_p}$  minimizes the residual and is obtained by solving the least squares problem 
\[ \min_{\bfy \in \R^{mn_p}} \| \bfr_m(\sigma) \|_2 =   \| \|\bfd\|_2 \bfe_1 - \bar{\bfH}_m(\sigma; \bfT_m)\bfy \|_2.\]
The details of this algorithm are given in~\cite[Algorithm 2]{bakhos2017multipreconditioned}.

\paragraph{Computational Cost} Assume that the cost of matvecs involving $\bfA_1$ and $\bfA_2$ is $T_{\rm matvec}$ and the cost of applying the preconditioner is $T_{\rm prec}$. If $k_m$ iterations of the MPGMRES-Sh are performed, then the total cost is 
\[ T_{\rm cost} =  k_mn_p(T_{\rm matvec} + T_{\rm prec}) + \mc{O}( N_\sigma n k_m^2 n_p^2) + \mc{O}(N_\sigma k_m^3n_p^3) \> \text{flops}. \]
As we shall see in the numerical experiments, the number of iterations $k_m$ low (typically $\leq 50$). Note that the cost of building the basis is independent of the number of shifted systems, whereas the cost of the projected system depends on the number of shifted systems $N_\sigma$.

\subsection{Applying the covariance operator}\label{ssec:applying} We now show how to compute the application of the covariance operator $\bfC_{\alpha}$ to a vector $\bff$. Let the matrices $\bfK$ and $\bfM$ be as defined in Section~\ref{ssec:disc}. 

Let $\alpha \in \R_+$ be a positive number. If $\alpha \in \mathbb{N}$, then it is straightforward to compute an application of the covariance operator which requires $\alpha$ sequential applications of solves involving $\bfK$. Henceforth, we assume that $\alpha \in \R_+ \setminus \mathbb{N}$. Let $r$ be the integer part of $\alpha$ and let $s \in (0,1)$ be the fractional part. We can apply the covariance operator to a vector {$\bff$, this} can be computed as
\begin{equation}\label{eqn:covapply} \bfC_{\alpha}\bff = \left(\sum_{j=-M_-}^{M_+} \frac{w_j}{z_j}   \left(  \bfM + z_j^{-1}\bfK \right)^{-1}\bfM\right) \left[\bfK^{-1} \bfM\right]^{r}\bff.\end{equation}
 In writing this expression, we have replaced the weights with $w_jz_j$. This form is completely equivalent to~\eqref{eqn:discov} but makes it easier to  precondition. To turn this expression into an efficient algorithm, we perform the following steps. In the offline phase, we determine the quadrature weights $\{w_j\}_{-M_-}^{M_+}$ and nodes $\{z_j\}_{-M_-}^{M_+}$, preconditioner shifts $\{\tau_j\}_{j=1}^{n_p}$. Furthermore, we compute factorizations of $\bfK$ and the preconditioners $\{ \bfP_j = \bfM + \tau_j\bfK\}_{j=1}^{n_p}$. In the online stage, we compute the action of $\bfK^{-1}\bfM$ on $\bff$ $r$ times to obtain $\bfc = [\bfK^{-1}\bfM]^r\bff$. We then apply MPGMRES-Sh with $\bfA_1 = \bfM$, $\bfA_2 = \bfK$, $\bfd = \bfM \bfc$ and shifts $\sigma_j = 1/z_j$ for $-M_{-} \leq j \leq  M_+ $ to obtain the solutions $\{\bfx_j\}_{-M_-}^{M_+}$.   Finally, we compute $\bfC_{\alpha}\bff = \sum_{j=-M_-}^{M_+} w_j \bfx_j$. The details are given in Algorithm~\ref{alg:fractional}.

When the integer part $r$ is large, the number of solves with $\bfK$ can become expensive; however, in practice, a value of $r > 3$ is typically not used since the resulting field is extremely smooth. If an application requires the use of a large value of $r$ techniques from~\cite[Section 4.1]{higham2008functions} may be used.

\paragraph{Alternative approaches} In the expression for the application of the covariance defined in~\eqref{eqn:covapply}, we used a slight reformulation for the fractional part, which is obtained by factoring out $z_j$ in each summand.  We also investigated two other formulations for computing $\bfD\bfb$. The first is perhaps the more natural version 
\[ \bfD\bfb =  \sum_{j=-M_-}^{M_+} w_j   \left(  z_j \bfM + \bfK \right)^{-1}\bfM\bfb. \]
While mathematically equivalent to the approach used in~\eqref{eqn:covapply}, we found that it was easier to precondition the formulation in~\eqref{eqn:covapply}. The second approach is to decompose 
\[ \bfD\bfb =  \sum_{j=-M_-}^{0} w_j(z_j \bfM + \bfK )^{-1}\bfM\bfb  + \sum_{j=1}^{M_+} \frac{w_j}{z_j}   \left(   \bfM + z_j^{-1}\bfK \right)^{-1}\bfM\bfb. \]
The advantage of this approach is that in this instance both the sets of weights, $w_j$ and $w_j/z_j$, are at most one in their respective intervals $-M_-\leq j \leq 0$ and $1 \leq j \leq M_+$; similarly the nodes $z_j$ and $z_j^{-1}$ are also at most $1$ in their respective intervals. However, the downside is that for apply the fractional operator, we need two solves with MPGMRES-Sh. For this reason, we prefer the formulation in~\eqref{eqn:covapply}.

In the offline stage Algorithm~\ref{alg:fractional}, we factorize the preconditioners. However, this may not be feasible for very large-scale problems, especially in three spatial dimensions. In such cases, we can apply the preconditioners using iterative solvers. See~\cite[Section 4.2.2]{bakhos2017multipreconditioned}.

\begin{algorithm}[!ht]
\begin{algorithmic}[1]
\REQUIRE matrices $\bfK$ and $\bfM$, vector $\bff$, scalar $\alpha \notin \mathbb{N}$, preconditioner shifts $\{\tau_{j}\}_{j=1}^{n_p}$
\STATE \COMMENT{Offline stage: Precomputation}
\STATE Split $\alpha = r + s$ where $s \in (0,1)$ and $r \in \mathbb{N}$
\STATE Determine quadrature weights $w_j$ and nodes $z_j$ for $M_- \leq j \leq M$ 
\STATE Factorize  $\bfP_j = \bfM + \tau_j \bfK$ for  $1 \leq j \leq n_p$  and $\bfK$

\STATE \COMMENT{Online stage: Solve} 
\STATE  Apply the integer part: $\bfc = \left[ \bfK^{-1} \bfM\right]^{r}\bff$  
\STATE Solve the sequence $(\bfM + z_j^{-1} \bfK) \bfx_j = \bfM\bfc$ and compute $\bfu = \sum_{j=M_-}^{M_+} w_j \bfx_j$. 
\RETURN Return $\bfu = \bfC_{\alpha}\bff$
\end{algorithmic}
\caption{Applying the covariance operator $\bfC_{\alpha}\bff$}
\label{alg:fractional}
\end{algorithm}

\subsubsection{Generating samples from the SPDE} We can adapt the efficient approach in Algorithm~\ref{alg:fractional} to generate samples from the SPDE
\begin{equation}\label{eqn:spdeomega} (\kappa^2(\bfx)  - \nabla\cdot (\bfTheta (\bfx) \nabla ))^{(\nu + d/2)/2} u(\bfx) = W(\bfx) \qquad \bfx \in \Omega,\end{equation}
with zero Neumann boundary conditions. We follow the approach in~\cite{bolin2018weak} to solve the SPDE. Let $\bfL\bfL^T = \bfM$ be the Cholesky factorization of $\bfM$ and let $\bfw \sim \mathcal{N}(\bfzero,\bfI)$. For $0 < \alpha/2 < 1$, we can compute a solution to the SPDE as 
\[\bfu := \frac{2\zeta\sin(\pi\alpha/2)}{\pi} \sum_{j=-M_-}^{M_+} e^{(1-\alpha/2)\zeta j}  ( \bfK + e^{j\zeta} \bfM)^{-1}\bfL\bfw,\] 
where $\zeta = 1/\log(1/h)$,  and the number of shifted systems  $N_\sigma = M_- + M_+ + 1$, where
\[ M_+ = \lceil {\pi^2}/{(2\alpha\zeta^2)} \rceil, \quad \text{and} \quad 
M_- = \lceil {\pi^2}/(2(2-\alpha)\zeta^2)\rceil.\] 
We can readily modify Algorithm~\ref{alg:fractional} to efficiently solve this problem. In particular, we can leverage the use of the MPGMRES-Sh approach to accelerate the solution to the SPDE. If we instead use zero Dirichlet boundary conditions for $u$, the error analysis in~\cite[Theorem 2.1]{bolin2018weak} applies directly.  

\subsection{Karhunen-Lo\`eve expansions}\label{ssec:kle} In this section, we discuss an approach generate samples from the Gaussian measure $\mc{N}(\mu, {\mc{C}})$ on a Hilbert space $\mc{H}$, where $\mu \in \mc{H}$ is a the mean function and  $\mc{C}$ is a self-adjoint, positive semidefinite covariance operator.  This approach is based on the Karhunen-Lo\`eve (KL) expansion of the stochastic process. In this approach we consider an orthonormal set of eigenpairs $(\lambda_j, \psi_j)$ for $1 \leq j \leq \infty$ where the eigenvalues are arranged in decreasing order as $\lambda_1 \geq \lambda_2 \geq \cdots$. Furthermore, let $\{\xi_j\}_{j=1}^\infty$ be an independent and identically distributed (i.i.d.) sequence of random variables with $\xi_1 \sim \mc{N}(0,1)$. Then by~\cite[Theorem 6.19]{stuart2010inverse}, the sample $u$ defined by the KL expansion
\[ u := \mu + \sum_{k=1}^\infty \sqrt{\lambda_j} \xi_j \psi_j, \]
is distributed according to $\mc{N}(\mu,{\mc{C}})$. In practice, for certain applications, the eigenvalues exhibit rapid decay and, therefore, the expansion can be approximated by the truncated representation $u_N := \mu + \sum_{j=1}^N \sqrt{\lambda_j} \xi_j \psi_j$. In the context of Bayesian inverse problems, rather than estimating the field, one can estimate the coefficients $\{ \xi_j\}_{j=1}^N$of the KL expansion. This can be computationally beneficial because of the input dimensionality reduction.

\paragraph{Galerkin approximation} In this discussion, we will derive the expressions for the truncated KL expansion for the Gaussian measure $\mc{N}(0, {\mc{C}})$. We will then discuss implementation details for {computing} the truncated KL expansion for the Gaussian measure $\mc{N}(0,\mc{C}_{\alpha})$. Consider a finite dimensional subspace $\mc{V}_h \subset L^2(\Omega)$ with $N_h = \text{dim}(\mc{V}_h)$ and let $\{\phi_j \}_{j=1}^{N_h}$ be a basis for the subspace. The eigenvalue problem $\mc{C} \psi = \lambda \psi$ is then replaced by the finite dimensional eigenvalue problem $\mc{C}_h \psi_h = \lambda_h \psi_h$, where $\psi_h \in \mc{V}_h$ and $\mc{C}_h: \mc{V}_h \rightarrow \mc{V}_h$. We seek a solution $\psi_h = \sum_{j=1}^{N_h} \psi_{h,j} \phi_j$, by considering the Galerkin projection
\[ \sum_{j=1}^{N_h}\psi_{h,j} (\phi_i, {\mc{C}}_h \phi_j)_{L^2(\Omega)} =  \sum_{j=1}^{N_h}\psi_{h,j} (\phi_i,  \phi_j)_{L^2(\Omega)} \qquad 1 \leq i \leq N_h. \]
Following the discussion in Section~\ref{ssec:disc}, we can write $(\phi_i,  \phi_j)_{L^2(\Omega)} = \langle \bfe_i,\bfe_j\rangle_{\bfM}$ and $(\phi_i, {\mc{C}}_h \phi_j)_{L^2(\Omega)} = \langle \bfe_i, \bfC\bfe_j\rangle_{\bfM}$, where $\bfC$ is the discretized covariance matrix. Therefore, we now have the generalized Hermitian eigenvalue problem (GHEP)
\begin{equation}\label{eqn:klghep} \bfM \bfC \bfpsi_h = \lambda_h \bfM \bfpsi_h, \end{equation}
where $\bfpsi_h = \bmat{\psi_{h,1} & \dots & \psi_{h,N_h}}^T\in \R^{N_h} $ is the generalized eigenvector. Let $\lambda_{h,1} \geq \dots  \geq \lambda_{h,N}$ be the generalized eigenvalues and let $\bfpsi_{h,j}$ be the corresponding eigenvectors. Then a sample from the process can be generated as 
\begin{equation}\label{eqn:disckle}
\bfu_{h,N} = \sum_{j=1}^N\sqrt{\lambda_{h,j}} \xi_j \bfpsi_{h,j}.
\end{equation}

\paragraph{Efficient implementation} We want to compute the truncated KL for the measure $\mc{N}(0,\mc{C}_{\alpha}).$ The GHEP~\eqref{eqn:klghep} can be solved using any standard Krylov-based eigensolver (e.g., Lanczos) that does not form the matrix $\bfC_{\alpha}$ or $\bfM\bfC_{\alpha}$ explicitly but instead relies on forming the matvecs involving the matrices. The matvecs involving $\bfC_{\alpha}$ can be efficiently computed using the MPGMRES-Sh approach as described in Algorithm~\ref{alg:fractional}.

An alternative approach to computing the truncated KL expansion is to work with the eigenpairs of the operator $\mc{A}$, which is attractive at first glance, since it avoids the complications with computing the fractional part. However, the approach based on the covariance operator $\mc{C}_{\alpha}$ directly is advantageous for two reasons. First, to compute the truncated KL expansion, we need to target the eigenpairs of $\mc{C}_{\alpha}$ corresponding to the largest eigenvalues, which corresponds to the eigenpairs of $\mc{A}$ corresponding to the smallest eigenvalues, which is typically much harder. Second, for $\alpha > 1$ the operator $\mc{A}^{-\alpha}$ acts as a spectral transformation that enhances the eigenvalue gaps.

\section{Application to Bayesian Inverse Problems}
In this section, we show how to leverage the efficient representation of the covariance operator. In Section~\ref{ssec:bayes}, {we review} the necessary background for Bayesian inverse problems. In Section~\ref{ssec:discpost}, we derive the discretized posterior distribution that uses the covariance matrix in~\eqref{eqn:discov} as the  prior covariance and in Section~\ref{ssec:genhybr}, we show how to adapt Generalized Hybrid Iterative Methods (GenHyBR) appropriately {to efficiently} compute the MAP estimate. 
\subsection{Bayesian Inverse Problems}\label{ssec:bayes}
Let $m \in L^2(\Omega)$ be the inversion parameter that we wish to recover from the measurements $\bfy \in \R^{N_y}$, which are related through the measurement equation
\begin{equation}
    \bfy = \mc{F}(m) + \bfeta,
\end{equation}
where $\mc{F}:L^2(\Omega) \rightarrow \R^{N_y}$ is the parameter-to-observable map or the measurement operator, and $\bfeta$ is the measurement noise which is assumed to be Gaussian with zero mean and covariance $\bfGamma_{\rm noise}$, which we write as $\mc{N}(\bfzero,\bfGamma_{\rm noise})$. Therefore, the likelihood probability density function takes the form
\[ \pi(\bfy| m) \propto \exp\left(-\frac12 \|\bfy - \mc{F}(m) \|_{\bfGamma_{\rm noise}^{-1}}^2\right). \]

We assume that $m$ is endowed with a Gaussian prior with mean function $m_{\rm pr}$ and covariance operator $\lambda_{\mc{C}}^{-2} \mc{C}_{\rm pr} $, i.e., $m \sim \mc{N}(m_{\rm pr},\lambda_{\mc{C}}^{-2}\mc{C}_{\rm pr})$. Here $\lambda_{\mc{C}}^2$ acts a regularization parameter that may be known, or can be estimated as part of the inversion process (this is what we do in this paper). For the operator $\mc{C}_{\rm pr}$ we take it to be $\mc{C}_{\alpha}$ defined in Section~\ref{ssec:prior}.

By the choice of the likelihood and the prior distribution, using the infinite dimensional Bayes formula, the Radon-Nikodym derivative of the posterior probability measure $\mu_{\rm post}$ with respect to the prior measure $\mu_{\rm pr} = \mc{N}(m_{\rm pr}, \lambda_{\mc{C}}^{-2}\mc{C}_{\rm pr})$, takes the form
\[ \frac{d\mu_{\rm post}}{d\mu_{\rm pr}}  =\frac{1}{Z} \exp \left( -\frac12 \|\bfy - \mc{F}(m) \|_{\bfGamma_{\rm noise}^{-1}}^2  \right),\]
where  $Z = \int_{}\pi(\bfy| m) d\mu_{\rm pr}$ is a normalization constant. 
The maximum a posteriori (MAP) estimator maximizes the posterior distribution and can be obtained by the solution to the optimization problem
\[ \min_{m \in L^2(\Omega)}\frac12 \|\bfy - \mc{F}(m) \|_{\bfGamma_{\rm noise}^{-1}}^2 + \frac{\lambda_{\mc{C}}^2}{2} \| m - m_{\rm pr}\|_{\mc{C}_{\rm pr}^{-1}}^2 . \]
If $\mc{F}$ is a linear operator, then the posterior distribution is Gaussian; we exploit this fact in this paper.

\subsection{Discretized posterior distribution}\label{ssec:discpost}

The covariance matrix $\bfC_{\rm pr}$ after discretization takes the same form as in~\eqref{eqn:discov}.
It is easy to verify that the resulting matrix $\bfC_{\rm pr}$ is self-adjoint; i.e., $\bfC_{\rm pr}^* = \bfC_{\rm pr}$ and $\bfC_{\rm pr}\bfM^{-1}$ is symmetric with respect to the standard Euclidean inner product and positive definite. This gives us a discretized representation of the prior distribution $\bfm \sim \mc{N}(\bfm_{\rm pr},  \lambda_{\mc{C}}^{-2}\bfC_{\rm pr}\bfM^{-1})$, i.e.,
\[ \pi_{\rm pr}(\bfm) \propto \exp\left( -\frac{\lambda_{\mc{C}}^2}{2} \langle \bfm-\bfm_{\rm pr}, \bfC_{\rm pr}^{-1} (\bfm-\bfm_{\rm pr}) \rangle_{\bfM} \right). \]

For completeness, we give expressions for the posterior distribution, and the MAP estimate, for the discretized problem. This uses the components developed in Section~\ref{ssec:disc}. The posterior probability density function takes the form 
\[ \pi_{\rm post}(\bfm|\bfy) \propto \exp\left( -\frac12 \|\bff(\bfm) - \bfy\|_{\bfGamma_{\rm noise}^{-1}}^2 -\frac{\lambda^2_{\mc{C}}}{2} \|\bfm-\bfm_{\rm pr} \|_{\bfM\bfC_{\rm pr}^{-1}}^2 \right). \]
The MAP estimator can be obtained by solving the optimization problem
\[ \min_{\bfm \in \R^{N_h}}\frac12 \|\bff(\bfm) - \bfy\|_{\bfGamma_{\rm noise}^{-1}}^2 +\frac{\lambda^2_{\mc{C}}}{2} \|\bfm-\bfm_{\rm prior} \|_{\bfM\bfC_{\rm pr}^{-1}}^2.\]
If the forward operator is linear $\bff(\bfm) = \bfF \bfm$, then the posterior distribution is Gaussian with $\bfm|\bfy \sim \mc{N}(\bfm_{\rm post}, \bfC_{\rm post}\bfM^{-1})$, where

\begin{equation}\label{eqn:postmean} \begin{aligned} \bfC_{\rm post} := & \> (\bfF^{\sharp}\bfGamma_{\rm noise}^{-1} \bfF + \lambda_{\mc{C}}^2\bfC^{-1}_{\rm pr})^{-1} \\ \bfm_{\rm post} := & \> \bfC_{\rm post}  (\bfF^{\sharp}\bfGamma_{\rm noise}^{-1}\bfy + \lambda_{\mc{C}}^2\bfC_{\rm pr}^{-1}\bfm_{\rm pr}), \end{aligned}
\end{equation}
and $\bfF^\sharp : \R^{N_y} \rightarrow \R^{N_h}_{\bfM}$ is the adjoint operator of $\bfF$, and takes the form $\bfF^\sharp = \bfM^{-1}\bfF^T$. Note that $\bfC_{\rm post}$ is self-adjoint with respect to the $\langle\cdot, \cdot \rangle_{\bfM}$ inner product, so that $\bfC_{\rm post}\bfM^{-1}$ is symmetric positive definite.

\subsection{GenHyBR for computing the MAP estimate}\label{ssec:genhybr}
We now derive an alternate expression for the MAP estimate~\eqref{eqn:postmean}. Multiplying both sides of~\eqref{eqn:postmean} with $\bfC_{\rm post}^{-1}$, we get   
\[ (\bfF^\sharp\bfGamma_{\rm noise}^{-1} \bfF +\lambda_{\mc{C}}^2\bfC_{\rm pr}^{-1})\bfm_{\rm post} = \bfF^\sharp\bfGamma_{\rm noise}^{-1}\bfy + \lambda_{\mc{C}}^2\bfC_{\rm pr}^{-1}\bfm_{\rm pr}. \]
We use $\bfF^\sharp =\bfM^{-1} \bfF^T$, and the change of variables $\widehat{\bfm}_{\rm post} = \bfM\bfC_{\rm pr}^{-1} (\bfm_{\rm post} - \bfm_{\rm pr})$. Then to obtain the MAP estimate, we solve the linear system for $\widehat{\bfm}_{\rm post}$
\begin{equation}\label{eqn:mapeq} (\bfF^T\bfGamma_{\rm noise}^{-1} \bfF \bfC_{\rm pr}\bfM^{-1} + \lambda_{\mc{C}}^2\bfI ) \widehat{\bfm}_{\rm post} =  \bfF^T\bfGamma_{\rm noise}^{-1} (\bfy - \bfF\bfm_{\rm pr}), 
\end{equation}
and then compute $\bfm_{\rm post} = \bfm_{\rm pr} + \bfC_{\rm pr}\bfM^{-1}\widehat{\bfm}_{\rm post}$. This form of the MAP estimator will be the basis for our computationally efficient procedure. 
Note that $\bfC_{\rm pr}\bfM^{-1}$ is a symmetric positive definite matrix with respect to the Euclidean inner product; we denote this by $\bfQ = \bfC_{\rm pr}\bfM^{-1}$ for convenience. {The application} of $\bfQ$  to a vector can be performed efficiently, using a slight modification of Algorithm~\ref{alg:fractional}. Furthermore, note that $\bfF^T\bfGamma_{\rm noise}^{-1} \bfF \bfQ$ is symmetric with respect to the $\langle \cdot, \cdot \rangle_{\bfQ}$ inner product. 

We can reformulate the equations in such a way that we can use the generalized Golub-Kahan bidiagonalization (gen-GK) approach~\cite{chung2017generalized} for efficiently estimating the MAP estimate. We initialize the iterations with $\bfb = \bfy- \bfF\bfm_{\rm pr}$,  $\delta_1 = \|\bfb\|_{\bfGamma^{-1}_{\rm noise}}, \bfu_1 = \bfy_1/\delta_1$, and $\gamma_1\bfv_1 = \bfF^T\bfGamma^{-1}_{\rm noise}\bfb$. At step $k$ in this approach, 
\[\begin{aligned}\delta_{k+1}\bfu_{k+1} = & \> \bfF\bfQ\bfv_k - \gamma_k \bfu_k \\
\gamma_{k+1}\bfv_{k+1} = & \>\bfF^T\bfGamma_{\rm noise}^{-1}\bfu_{k+1} - \delta_{k+1}\bfv_k,
\end{aligned}
\]
where $\gamma_i,\delta_i \geq 0$ are chosen such that $\|\bfu_i\|_{\bfGamma_{\rm noise}^{-1}} =\|\bfv_i\|_{\bfQ} = 1 $ for $1 \leq i \leq k+1$. These iterates can be collected to form the matrices $\bfU_{k+1} := \bmat{\bfu_1 & \dots & \bfu_{k+1}}$, $\bfV_k := \bmat{\bfv_1 & \dots & \bfv_k}$, and the bidiagonal matrix
\[ \bfB_k := \bmat{\gamma_1 \\ \delta_1 & \gamma_2 \\ & \delta_2 & \ddots \\ & & \ddots & \gamma_k \\ & & &  \delta_{k+1} } \in \R^{(k+1)\times k}. \]
This can be rearranged to obtain the relations, which are typically accurate up to machine precision
\[ \begin{aligned} \bfF\bfQ\bfV_k = & \> \bfU_{k+1} \bfB_k \\ \bfF^T\bfGamma_{\rm noise}^{-1}\bfU_{k+1} = & \> \bfV_k \bfB_k  + \gamma_{k+1}\bfv_{k+1}\bfe_{k+1}^T,\end{aligned}\]
and the orthogonality relations $\bfU_{k+1}^T\bfGamma_{\rm noise}^{-1}\bfU_{k+1} = \bfI_{k+1}$ and $\bfV_k^T\bfQ\bfV_k = \bfI_k$. The columns of $\bfV_k$ form a basis for the Krylov subspace $\mc{K}_k(\bfF^T\bfGamma_{\rm noise}^{-1}\bfF\bfQ, \bfF^T\bfGamma_{\rm noise}^{-1}\bfb)$, where the Krylov subspace is defined as $\mc{K}_k(\bfJ,\bfd) := \text{span}\{\bfd, \bfJ\bfd, \dots,\bfJ^{k-1}\bfd\}.$ To obtain the approximate solution, we search for linear combinations of the columns $\bfV_k$, i.e., $\widehat\bfm_{k,\lambda} = \bfV_k \bfz_{k,\lambda}$ that solve the optimization problem
\[ \min_{\bfz \in \R^k} \frac12 \|\bfB_{k+1}\bfz - \delta_1 \bfe_1\|_2^2 + \frac{\lambda_{\mc{C}}^2}{2} \|\bfz\|_2^2.  \]
To estimate the regularization parameter $\lambda_{\mc{C}}$, we minimize the projected Generalized Cross Validation (GCV). Other choices for estimating the regularization parameter are possible (e.g., discrepancy principle, unbiased predictive risk estimate), but we do not pursue them here. We refer the reader to~\cite{chung2017generalized} for additional details. To terminate the iterations we use a stopping criteria similar to~\cite[Section 5.5]{chung2008weighted}. Suppose the iterations terminate at step $k$, then the solution to the MAP estimate can be obtained by undoing the change of variables; that is, we compute the approximate solution $\bfm_{{\rm post}}^{(k)} :=  \bfm_{\rm pr} + \bfC_{\rm pr}\bfM^{-1}\widehat\bfm_{k,\lambda}$.

It is worth mentioning that each iteration of gen-GK requires two matrix-vector products (matvecs) with $\bfQ$, one with $\bfF$ and $\bfF^T$. There is an additional cost of $\mc{O}(k^3)$ for solving the projected least-squares problem and estimating the regularization parameter $\lambda_{\mc{C}}^2$, and an additional $\mc{O}(k(m+n))$ for orthogonalization. The additional cost of the extra matvec with $\bfQ$ is offset by the fact that we can efficiently estimate the regularization parameter during the iterative scheme. 
\subsection{Posterior Variance}\label{ssec:postvar} Methods to compute the approximate posterior variance have been developed in several references~\cite{flath2011fast,bui2013computational,villa2018hippylib,spantini2015optimal,saibaba2020efficient}. However, these are not directly applicable here and have to be suitably modified.  We follow the approach in~\cite{saibaba2020efficient}, that reutilizes the GenGK basis vectors computed during the computation of the MAP estimate, but we need to modify it in two different ways: first, since we are approximating infinite-dimensional quantities, we have to work with additional mass matrices, and, second, we no longer have ready access to the diagonals of $\bfQ$.

 The variance field corresponding to a covariance operator $\mc{C}$ can be {obtained} using the discussion in~\cite{bui2013computational}.  Once again suppose that we are given a basis $\{\phi_j \}_{j=1}^{N_h}$.   Given the discretized representation $\bfC$ of the covariance operator $\mc{C}$,  then the discretized covariance function 
\[ c_h(\bfx,\bfy) = [\bfPhi(\bfx)]^T \bfC\bfM^{-1} \bfPhi(\bfy), \]
where $\bfPhi(\bfx) = \bmat{\phi_1(\bfx), \dots, \phi_{N_h}(\bfx) }^T$ is a vector containing the finite element basis functions. To visualize the covariance function, we consider the discretized covariance function evaluated at the nodal points $\{\bfx_j\}_{j=1}^{N_h}$. Therefore, it is sufficient to consider the diagonals of the matrix $\bfC\bfM^{-1}$. For the posterior variance, therefore, we have to compute the diagonals of $\bfC_{\rm post}\bfM^{-1}$; that is, we need to compute the diagonals o 
\[ \bfC_{\rm post}\bfM^{-1} = (\bfF^{\sharp}\bfGamma_{\rm noise}^{-1} \bfF + \lambda_{\mc{C}}^2\bfC^{-1}_{\rm pr})^{-1} \bfM^{-1} = (\bfF^{T}\bfGamma_{\rm noise}^{-1} \bfF + \lambda_{\mc{C}}^2\bfQ^{-1})^{-1}. \] 
This can be efficiently estimated using the intermediate computations in the GenGK process~\cite{chung2018efficient,saibaba2020efficient}.  More precisely, we can approximate 
\[  \bfQ\bfF^T \bfGamma_{\rm noise}^{-1} \bfF\bfQ  \approx\bfQ \bfV_k \bfB_k^T \bfB_k  \bfV_k^T \bfQ = \bfZ_k \bfPhi_k \bfZ_k^T,\]
where $\bfB_k^T\bfB_k = \bfW_k\bfPhi_k\bfW_k^T$ is the eigendecomposition of $\bfB_k^T\bfB_k$ and $\bfZ_k = \bfQ\bfV_k\bfW_k$. Then we have the following approximation to the posterior covariance 
\[ \bfC_{\rm post} \approx \widehat{\bfC}_{\rm post} :=  \bfQ( \bfZ_k\bfPhi_k\bfZ_k^T + \lambda_{\mc{C}}^2\bfQ)^{-1} \bfQ \bfM. \]
Applying the Woodbury identity, we get 
\[ \widehat{\bfC}_{\rm post} =  \lambda_{\mc{C}}^{-2} \bfQ \bfM - \bfZ_k\bfDelta_k\bfZ_k^T \bfM ,\]
where $\bfDelta_k \in \R^{k\times k}$ is a diagonal matrix with the $i$th diagonal $\lambda_{\mc{C}}^{-2}\phi_i/(\phi_i + \lambda_{\mc{C}}^{2})$ for $1 \leq i \leq k$. Finally, we can get the approximate posterior variance from the diagonals of 
\[\widehat{\bfC}_{\rm post}\bfM^{-1}  =\lambda_{\mc{C}}^{-2} \bfQ  - \bfZ_k\bfDelta_k\bfZ_k^T . \]
Since $\bfZ_k\bfDelta_k\bfZ_k^T$ is low-rank, its diagonals can be easily computed. Finally, since in floating point arithmetic, the vectors $\bfV_k$ lose orthogonality, we use full reorthogonalization which adds additional cost; see~\cite{chung2018efficient,saibaba2020efficient}.  To estimate the diagonals of $\bfQ$ we follow the Diag++ approach in~\cite[Algorithm 1]{baston2022stochastic}. This is a Monte-Carlo based approach that only uses matvecs with $\bfQ$ (see Algorithm~\ref{alg:fractional}) to estimate the diagonals.

\section{Numerical Experiments}
We present a suite of numerical experiments that demonstrate the performance of the proposed methods. The timing results are reported on a Mac Mini (M1, 2020) with $16$ GB memory and running macOS Big Sur 11.2.3 and MATLAB 2021A. 
\subsection{Application of the covariance operator}
We perform some numerical experiments demonstrating the efficiency of the MPGMRES-Sh solver in computing the action of the covariance operators. We take the domain to be $\Omega = (0,1)^2$ and the number of grid points to vary from $33\times 33$ to  $513 \times 513$. We take $\mc{C}_{\alpha} = \mc{A}^{-s}$ where $\mc{A} = \kappa^2  - \Delta$, $\alpha = s \in (0,1)$. We also take $\kappa^2 = 100$ and the Laplacian has zero Neumann boundary conditions.  Although this operator is not trace-class, the goal here is merely to study the performance of MPGMRES-Sh and the results are applicable to other values of $\alpha > d/2$. 

Since the shifts $\sigma_j = 1/z_j$ are on the positive real axis, we use three preconditioners $n_p = 3$ with $\tau \in \{ 10^{-8}, 10^{-4}, 10^{-2}\}$. The preconditioners $\bfP_j$ for $j=1,\dots,n_p$ are factorized and stored ahead of the MPGMRES-Sh iterations. We stopped the MPGMRES-Sh iterations when the relative residuals of each shifted system were smaller than $10^{-8}$. We investigated with different numbers of shifts and with different shifts, but found this setup to strike a balance between the number of iterations and the cost per iteration.

\begin{figure}[!ht]\centering
\includegraphics[scale=0.25]{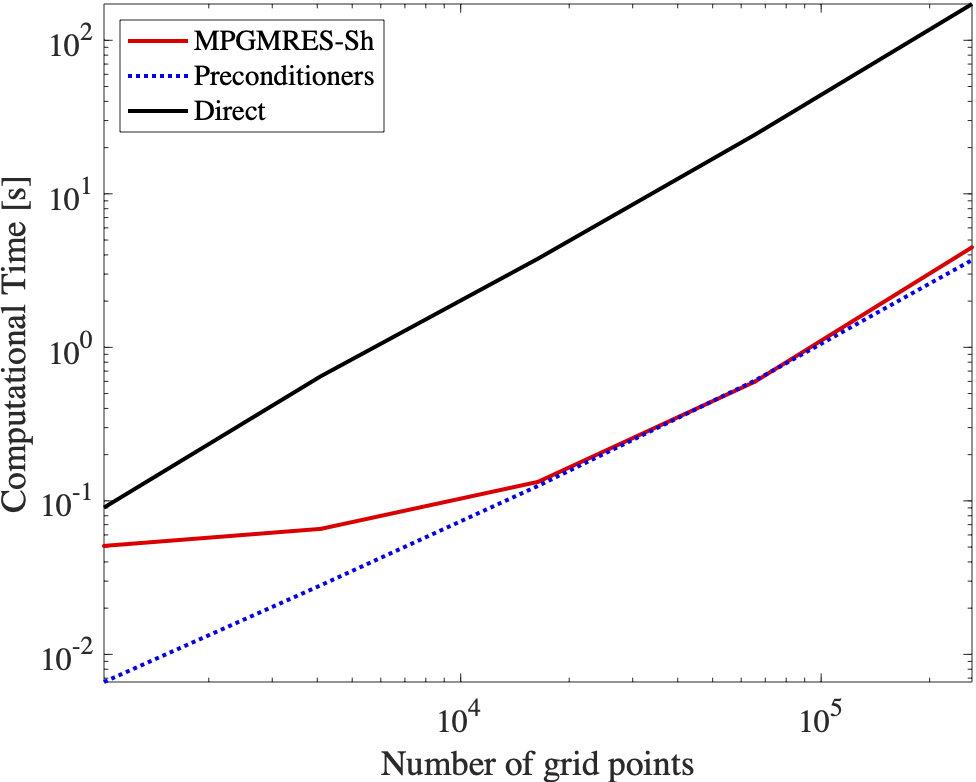}
\includegraphics[scale=0.25]{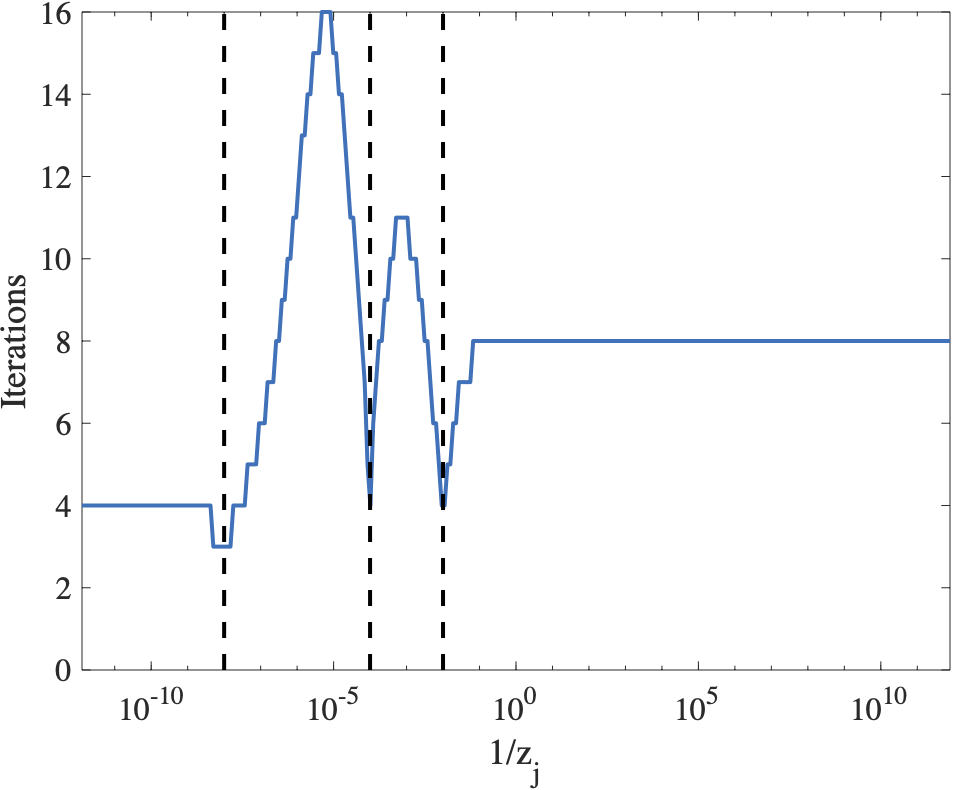}
\caption{(left) Timing results with increasing discretization for computing $\bfC_\alpha \bff$ where $\alpha = 0.5$. `Direct' refers to solving the sequence of shifted linear systems individually using a direct solver, `Preconditioner' refers to the time required to construct the preconditioners  $\{\bfP_j\}_{j=1}^{n_p}$. (right) The number of iterations at convergence taken by MPGMRES-Sh for each shift; note that a single basis was used across all the systems. The vertical black lines denotes the values of $\tau_j$ used for preconditioners. }
\label{fig:timing}
\end{figure}

\paragraph{Varying mesh discretization} In this experiment, we study the  effect of mesh refinement on the performance of the MPGMRES-Sh solver. We vary the grid sizes from $33\times 33$ to $513\times 513$. In the largest problem instance, the number of degrees of freedom was $263169$.  For the MPGMRES-Sh algorithm, we still used the same three shifts to determine the preconditioners. We see that even though the number of shifted systems to be solved increases with mesh refinement, the number of MPGMRES-Sh iterations only rises by a modest amount. The highest dimension of the  basis for the search space is $22\times 3 = 66$. In Figure~\ref{fig:timing}, we compare the timing cost of the MPGMRES-Sh with other approaches with mesh refinement. The `Direct' approach refers to solving the sequence of shifted systems using a direct solver. We also report the time taken in factorizing the preconditioners $\{\bfP_j\}_{j=1}^{n_p}$, labeled `Preconditioner', and the time taken for solving the shifted linear systems with MPGMRES-Sh is labeled `MPGMRES-Sh'. We see that the cost of MPGMRES-Sh is much smaller than the Direct method and is comparable to the cost of factorizing the preconditioners. For the largest problem size, there is a $\sim 40$X speedup of MPGMRES-Sh compared to Direct. The computational gains are more pronounced as the system size gets larger since the number of shifted systems grows significantly but the size of the basis only exhibits mild growth.  

\begin{table}[!ht]
    \centering
    \begin{tabular}{c|c|c}
    Grid size   & $N_\sigma$ &  Iters \\ \hline
     $33\times 33$ & $123$ & $9$ \\
     $65\times 65$ & $173$ & $10$ \\
     $129\times 129 $ & $235$ & $11$ \\
     $257\times 257$ & $305$ & $16$ \\ 
     $513 \times 513$ & $387$ & $22$     
      \end{tabular}
    \caption{Effect of the mesh discretization on the number of MPGMRES-Sh iterations. Also reported are the number of shifted systems to be solved. }
    \label{tab:mpgmres_h}
\end{table}

\begin{figure}[!ht]
    \centering
    \includegraphics[scale=0.35]{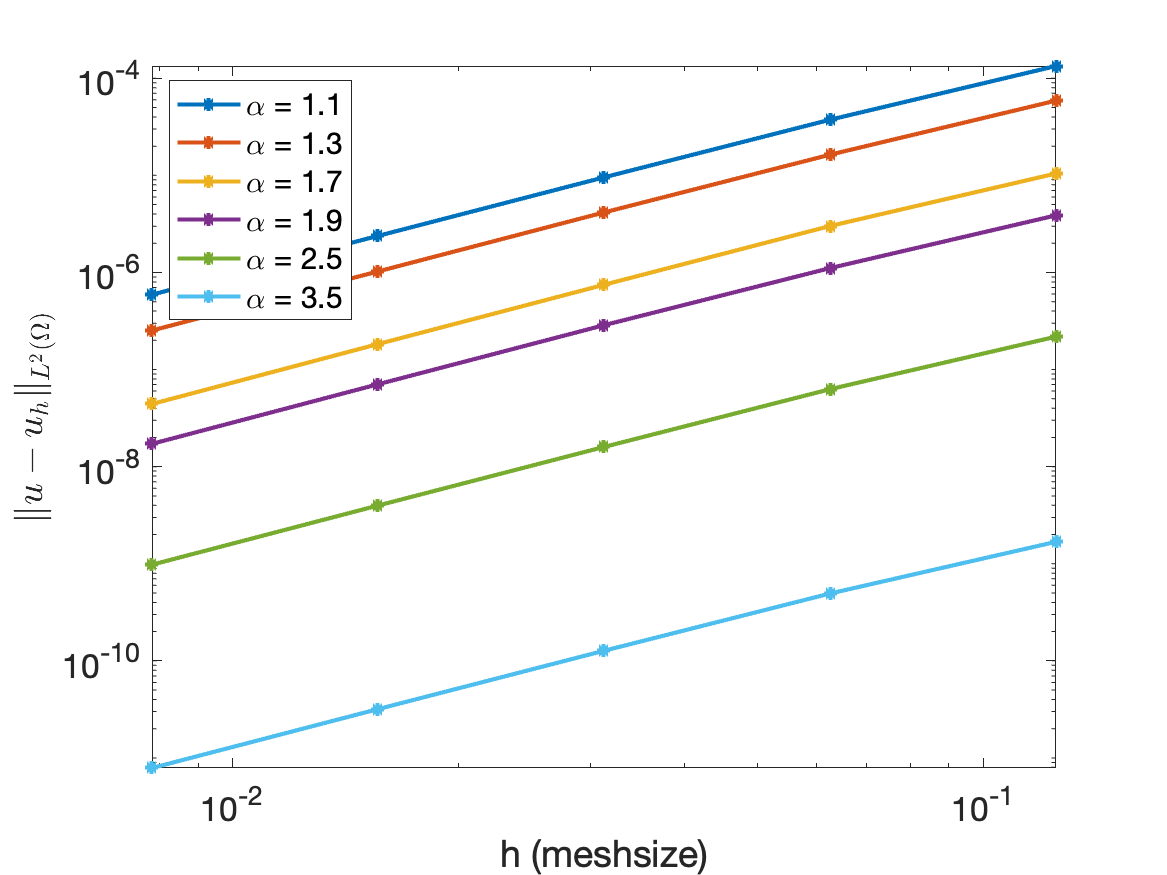}
    \caption{Accuracy of the covariance approximation for different values of $\alpha$ and mesh discretization $h$.}
    \label{fig:femaccuracy}
\end{figure}

Next, we compute the accuracy of the covariance approximation. We follow the method of manufactured solutions and let $f = \cos(2\pi x_1)\cos(2\pi x_2)$, so that $u = \mc{C}_\alpha f = (\kappa^2 + 8\pi^2)^{-\alpha}\cos(2\pi x_1)\cos(2\pi x_2)$. We compute the approximation using the formula in~\eqref{eqn:discov} computed using the approach in Algorithm~\ref{alg:fractional}. We then compute the absolute error in the $L^2(\Omega)$ norm and plot the error in Figure~\ref{fig:femaccuracy} for different values of $\alpha$ and different values of $h$, the mesh discretization parameter. We see that for all the values of $\alpha$, the error decreases with decreasing values of $h$ on the order of $h^2$. Furthermore, the absolute error also decreases with increasing values of $\alpha$. The error was analyzed for $0 < \alpha < 1$ in~\cite{bonito2015numerical} for the Dirichlet boundary case. {In order to apply such an analysis to our setting, this analysis needs to be first extended to Neumann boundary setting and then to the case $\alpha > 1$. We leave this as part of future work.} 

\paragraph{Effect of different values of $\alpha$}  In this experiment,  we take the grid size to be $257\times 257$ and vary the exponent $s$ from $0.1$ to $0.9$ in increments of $0.1$. The values of the number of shifted systems and the number of iterations taken by MPGMRES-Sh are given in Table~\ref{tab:mpgmres}. We see that the number of iterations is relatively small and it does not depend on the value of $s$. Furthermore, while the number of shifted systems can be large for each value of $s$, the number of iterations is relatively small. Since the same basis is used across all the shifted systems, the cost is amortized. Furthermore, the same basis can be used across multiple values of $s$ since the behavior is essentially independent of $s$. This can be very advantageous in a multi-query setting, where we need to apply the covariance operator to the same vector for multiple values of $s$.

\begin{table}[!ht]
    \centering
    \begin{tabular}{c|c|c|c|c|c|c|c|c|c}
     $s$      & $0.1$ & $0.2$ & $0.3$ & $0.4$ & $0.5$ & $0.6$ & $0.7$ & $0.8$ & $0.9$  \\ \hline
     $N_\sigma$    & $846$ & $476$ & $364$ & $318$ & $305$& $318$ & $364$ & $476$ & $846$ \\ 
     Iters & $16$ & $16$ & $16$ & $16$ & $16$ & $16$ & $16$ & $16$ & $16$
    \end{tabular}
    \caption{The number of shifted systems and the number of iterations required by MPGMRES-Sh. Three preconditioners are used corresponding to the shifts $\{10^{-8}, 10^{-4}, 10^{-2}\}$. }
    \label{tab:mpgmres}
\end{table}

\subsection{Samples from the SPDE}We provide numerical experiments for generating samples from the Gaussian process defined by the SPDE~\eqref{eqn:spde} with zero Neumann boundary conditions. We take the domain to be $\Omega =(0,1)^2$ and set the grid size to $129\times 129$. We solve the SPDE using the approach described in Section~\ref{ssec:applying} with MPGMRES-Sh used to accelerate the solutions of the shifted linear systems. We use the same preconditioners as in the previous experiment.  

\begin{figure}[!ht]
    \centering
    \includegraphics[scale=0.37,trim={2cm 0 2cm 0},clip]{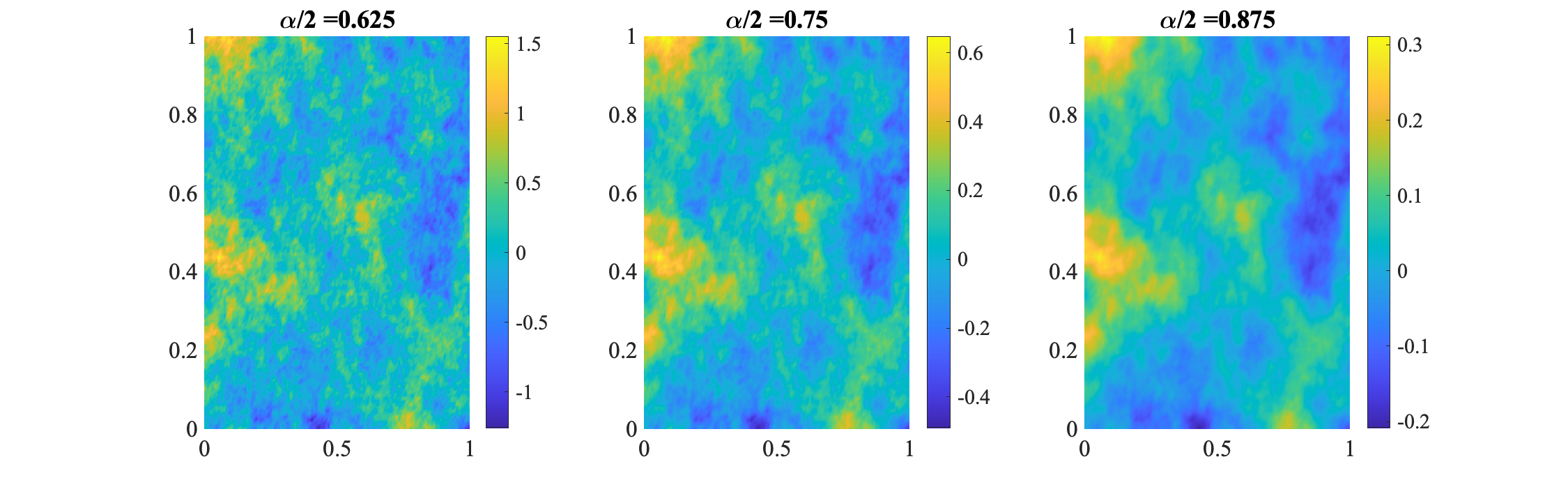}
    \caption{Samples from the SPDE~\eqref{eqn:spdeomega} with $\nu  \in \{\frac14, \frac24, \frac34\}$ and $\alpha = \nu + d/2$.}
    \label{fig:samples}
\end{figure}

We choose three different values of $\nu \in \{\frac14, \frac24, \frac34\}$ and take $\kappa^2 = 100$. Note that for these choices of $\nu$, we have  $1/2 < \alpha/2 < 1 $. 
 The samples are visualized in Figure~\ref{fig:samples}; note that we have used the same realization $\bfw$ across all the values of $\alpha$. It is readily seen that the for larger values of $\alpha$, the sample realizations appear smoother. The number of MPGMRES-Sh iterations for each value of $\alpha$ is the same and is $11$. 
 
\begin{figure}[!ht]
    \centering
    \includegraphics[scale=0.37,trim={2cm 0 2cm 0},clip]{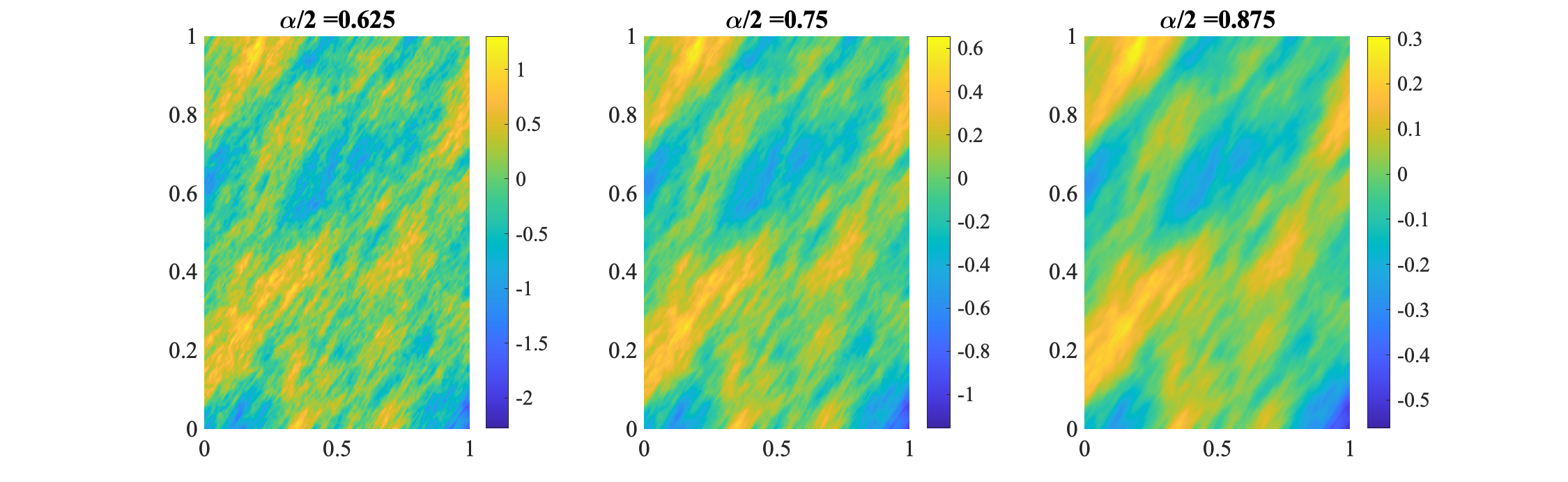}
    \caption{Samples from the anisotropic SPDE with $\nu  \in \{\frac14, \frac24, \frac34\}$ and $\alpha = \nu + d/2$. }
    \label{fig:spdeanisosamples}
\end{figure}

We also compute samples from the SPDE~\eqref{eqn:spdeomega}  with zero Neumann boundary conditions where, we take $\bfkappa^2(\bfx) = 100$ and 
 \begin{equation}\label{eqn:Theta}\bfTheta(\bfx) = \bmat{ \cos(\theta) & \sin(\theta) \\ - \sin(\theta) & \cos(\theta) } \bmat{\ell_1^2 \\ & \ell_2^2 }\bmat{ \cos(\theta) & -\sin(\theta) \\  \sin(\theta) & \cos(\theta) }  .   \end{equation}
 We take $\theta = \pi/4$, $\ell_1^2 = 10$, and $\ell_2 = 1$. Note that the corresponding random field is stationary and anisotropic. As before we pick $\nu  \in \{\frac14, \frac24, \frac34\}$ and the grid size is $129\times 129$. The samples are plotted in Figure~\ref{fig:spdeanisosamples}; note that we {have used} the same realization $\bfw$ across all the values of $\alpha$. The number of MPGMRES-Sh iterations taken for each value of $\alpha$ is the same and is $15$. 

\subsection{Computing the truncated KL expansion}
\begin{figure}[!ht]
    \centering
    \includegraphics[scale=0.3]{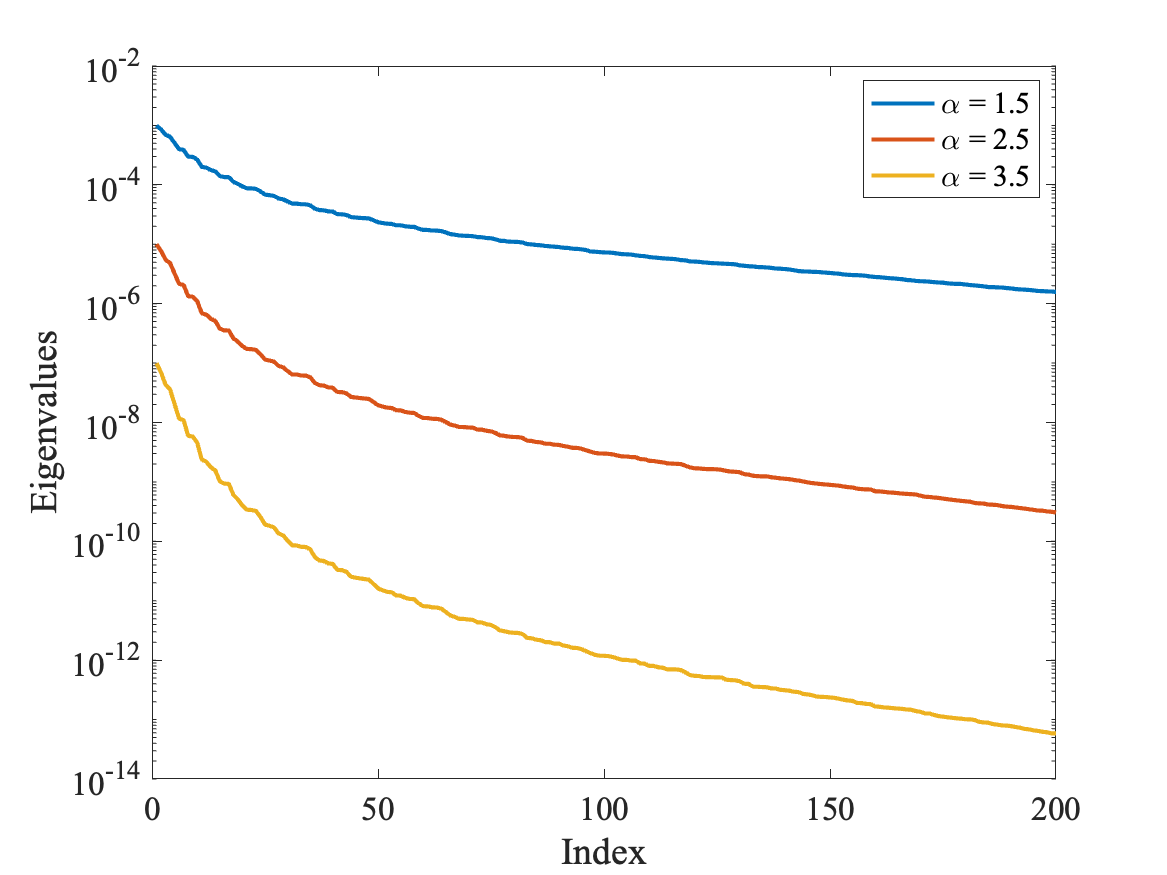}
    \includegraphics[scale=0.3]{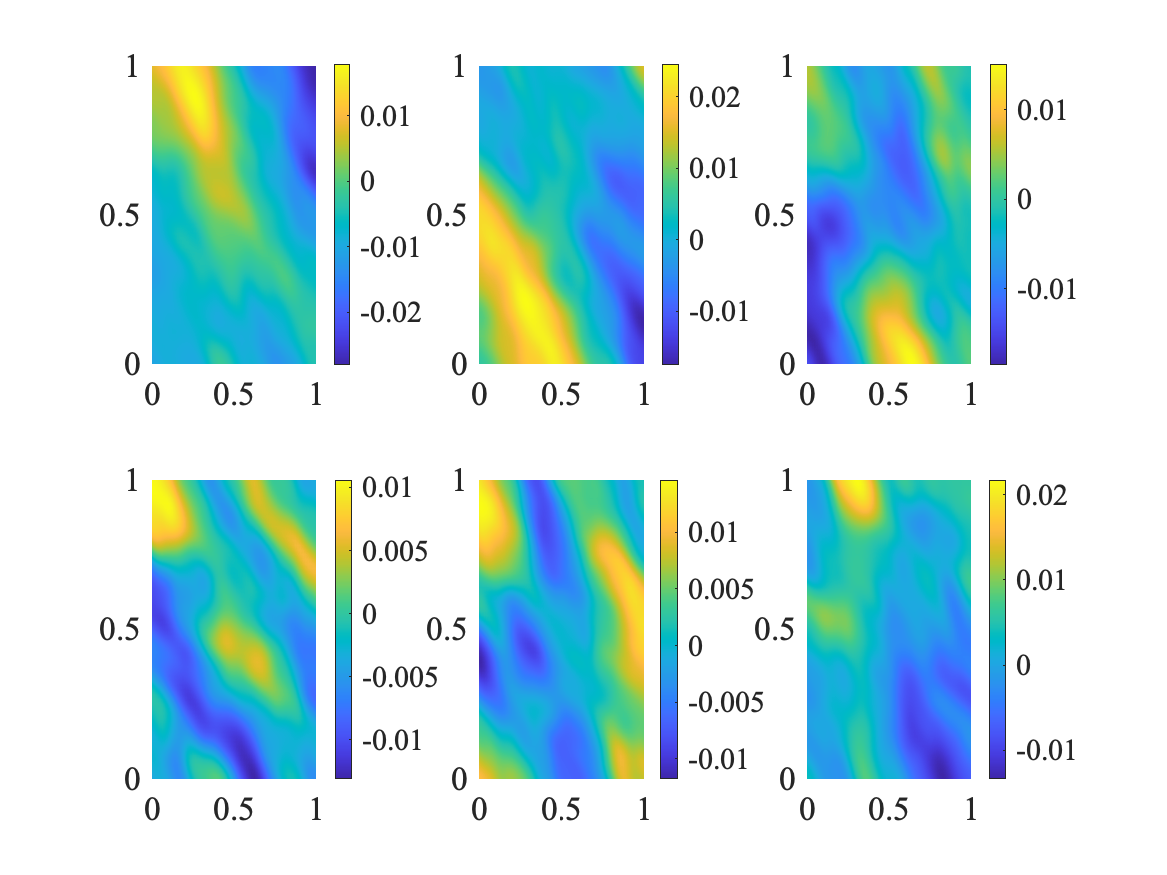}
    \caption{(left) Eigenvalues of the KL expansion for $\alpha \in \{\frac32,\frac52,\frac72 \}$, (right) 6 samples from the truncated KL expansion with $\alpha = \frac52$. }
    \label{fig:kle}
\end{figure}

We consider the anisotropic covariance operator~\eqref{eqn:covalpha} where $\kappa^2 = 80$ and $\bfTheta$ is defined in~\eqref{eqn:Theta}  with $\ell_1^2=4, \ell_2 = 1$ and $\theta = -\pi/4$. We compute the truncated KL expansion for three different values of $\alpha \in \{\frac32,\frac52,\frac72\}$. To compute the approximate eigenpairs, we used the two-pass randomized algorithm for GHEP~\cite[Algorithm 6]{saibaba2016randomized}. We computed $200$ eigenvalues and used an oversampling parameter of $20$. The eigenvalues are plotted in the left panel of Figure~\ref{fig:kle}; on the right panel, we plot $6$ different samples computed using the truncated KL expansion with $\alpha = 5/2$. 

\subsection{Tomography examples}\label{ssec:tomo}
\begin{figure}[!ht]
    \centering
    \includegraphics[scale=0.37,trim={2cm 0 2cm 0},clip]{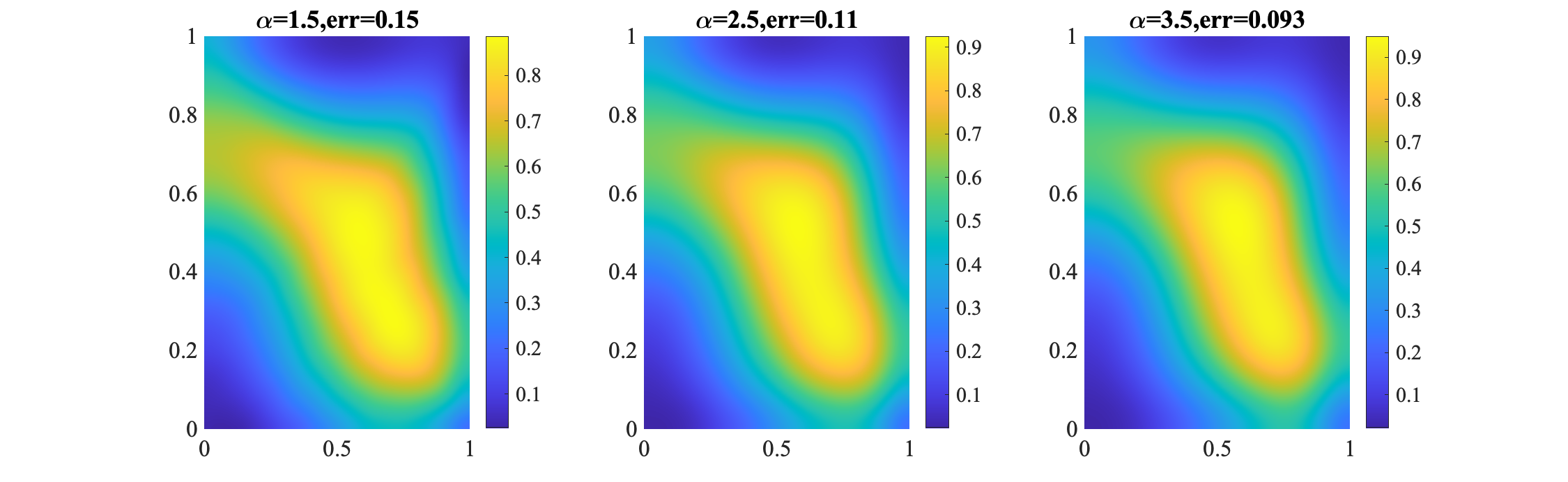}\\
    
    \includegraphics[scale=0.25]{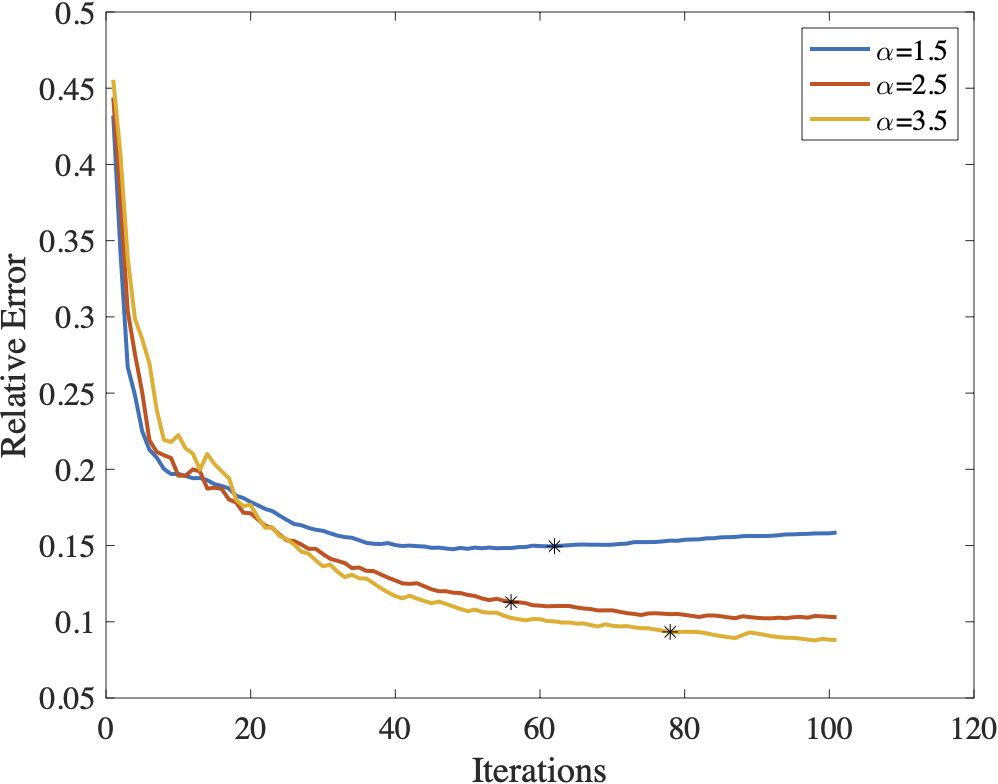}
    \caption{ Seismic tomography: (Top) Reconstructions for different values of $\alpha$; (bottom) iteration history of relative error. Black circles denote the stopping criterion.}
    \label{fig:seismic}
\end{figure}

\begin{figure}[!ht]
    \centering
    \includegraphics[scale=0.37,trim={2cm 0 2cm 0},clip]{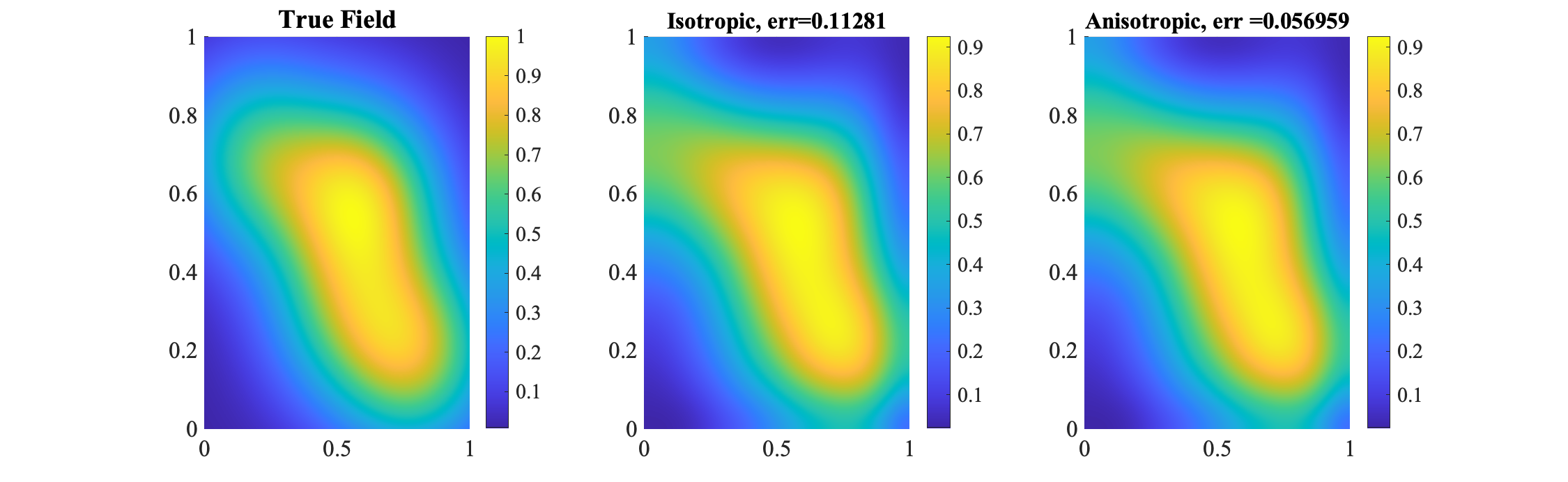}
    \caption{Seismic test problem: (left) true image, (center) reconstruction using Isotropic covariance, (right) reconstruction using Anisotropic covariance.}
    \label{fig:seismicaniso}
\end{figure}

We consider a {test problem} from the IRTools package~\cite{gazzola2019}. The number of grid points are $n_x \times n_y$. First, we use the \texttt{PRseismic} test problem with the number of sources $\lceil 0.4n_x\rceil$ and the number of receivers $\lceil 0.6 n_x \rceil$. In total, there were $4004$ measurements to which we add $2\%$ Gaussian noise. For the prior distribution with zero mean and the covariance matrix obtained from the discretization of $\mc{A}^{-\alpha}$ with $\mc{A} = \kappa^2  - \Delta$.  Next, we consider a test problem with real data from the Finnish Inverse Problems society~\cite{bubba2017tomographic}. In this instance of the test problem, the number of grid points is $128\times 128$ and the number of measurements is $14835$. We use the same prior distribution as before. A summary of the settings for the test problems is given in Table~\ref{tab:settings}.
\begin{table}[!ht]
\centering
\begin{tabular}{|c||c|c|}
\hline
 \textbf{Image}   &
   \textbf{Smooth}  &    \textbf{Cheese} \\
  
   \hline
   Application &  Seismic  & {X-ray }  \\
      \hline
   $m$  & $4004$   & $14835$ \\
    \hline
   $n_x \times n_y $ & $128 \times 128$   & $128 \times 128$ \\
    \hline
   $\kappa^2$  & $80$ & $300$ \\ \hline     
\end{tabular}
\caption{Problem settings for the tomography test problems in Section~\ref{ssec:tomo}.}
\label{tab:settings}
\end{table}
In Figure~\ref{fig:seismic}, we plot the results for the seismic test problem. In the top plot, we plot the reconstructions for three different values of $\alpha \in \{1.5,2.5,3.5\}$; since $d=2$, these correspond to $\nu= 1/2,3/2,5/2$. The relative error for each value of $\alpha$ is given in the title of each image. Since the images have smooth features, we see that the error has a slight decrease for increasing values of $\alpha$. In the bottom plot of the same figure, we plot the iteration history of the relative error computed using the GenHyBR method. We also highlight, in black, the iteration at which the algorithm has stopped. However, we plot the relative error until the maximum number of iterations $100$ to show that the relative error has stabilized and no semiconvergence is seen. 
In the next experiment, we consider the Seismic tomography problem but investigate the impact of an anisotropic covariance operator. We take the covariance operator~\eqref{eqn:covalpha}, with $\kappa^2 = 80$ and $\bfTheta$ defined in~\eqref{eqn:Theta}; we take $\ell_1 = 4, \ell_2 = 1$ and $\theta = -\pi/4$. For the isotropic covariance, we take $\kappa^2 = 80$ and $\bfTheta = \bfI$. For both the covariances, we take the value of $\alpha$ to be $2.5$. The true image, reconstruction with isotropic covariance and anisotropic covariances are displayed in Figure~\ref{fig:seismicaniso}. The relative reconstruction errors are given on the titles; for isotropic covariance GenHyBR took $56$ iterations and for anistropic covariance GenHyBR took $125$ iterations. It is clearly seen that by choosing the anisotropic covariance as the prior covariance, we can reduce the reconstruction error. This suggests that the anisotropic covariance operators may be beneficial under certain circumstances; however, to pick the correct parameters one has to use expert knowledge or they have to be estimated from data. 

\begin{figure}[!ht]
    \centering
    \includegraphics[scale=0.37,trim={2cm 0 2cm 0},clip]{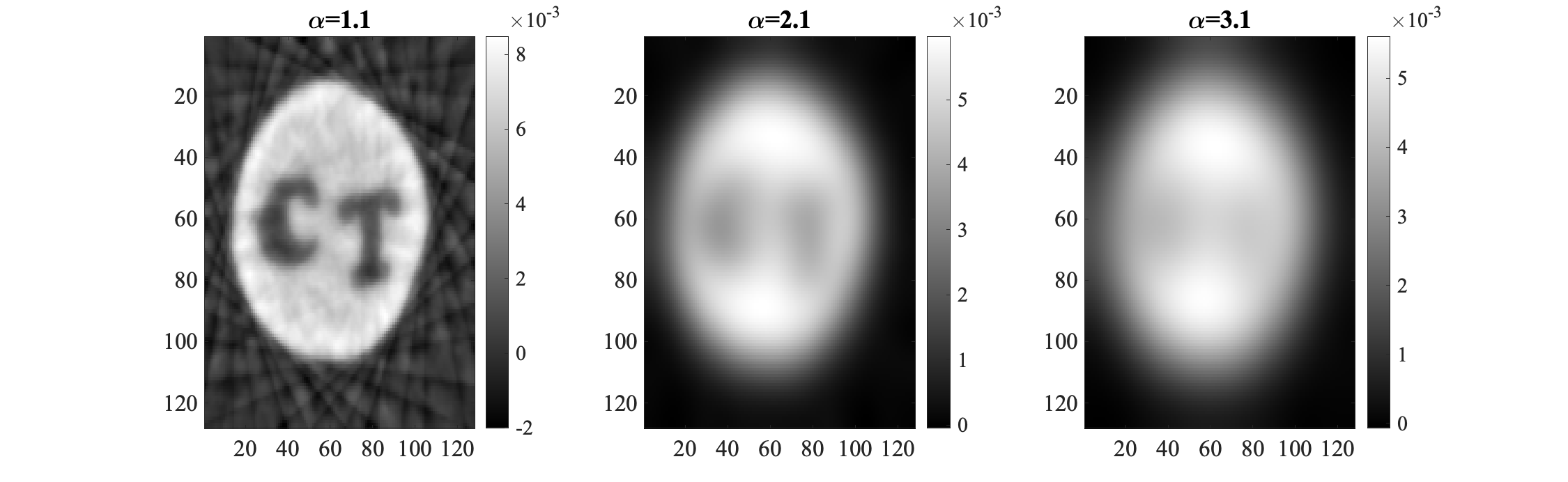}
    \caption{X-ray tomography: Reconstructions for different values of $\alpha$. }
    \label{fig:cheese}
\end{figure}

Finally, in Figure~\ref{fig:cheese}, we plot the result of the reconstructions corresponding to the real data problem using X-ray tomography. Since the true image is not available, we do not plot the relative error history. We choose $\alpha \in \{1.1,1.5,1.9\}$ since for higher values of $\alpha$, we get poor reconstructions. We have used the value of $\kappa^2 = 300$ and limit the maximum number of iterations to $150$; the number of iterations taken were $150$, $83$, and $75$ for $\alpha \in \{1.1,2.1,3.1\}$ respectively.

\subsection{PDE-based example} In this application, we consider a forward problem which is PDE-based. The underlying PDE is a 2D time-dependent heat equation
\begin{equation}
    \begin{aligned}
    \frac{\partial u}{\partial t} = & \Delta u  \qquad & \bfx \in \Omega, t \in {(0,T]}\\ 
    \bfn\cdot \nabla u = & 0 \qquad & \bfx \in \partial \Omega , t \in {(0,T]}\\
    u(t=0) = & \> m   \qquad & \bfx \in \Omega .
    \end{aligned}
\end{equation}
 Here $\Omega = (0,1)^2$, $T=0.01$, and $\bfn$ is the outward normal vector. The forward problem is discretized using Galerkin finite elements and Crank-Nicholson time differences; we use a grid size of $65\times 65$ with $100$ time steps. The inverse problem involves reconstructing the initial conditions from a discrete set of measurements at the final time point $t = T$. We add $2\%$ Gaussian noise to simulate measurement error. We use the covariance operator~\eqref{eqn:covalpha} {with} $\kappa^2 = 80$, $\bfTheta = \bfI$, and $\alpha = 2.5$.

\begin{figure}[!ht]
    \centering
    \includegraphics[scale=0.37,trim={2cm 0 2cm 0},clip]{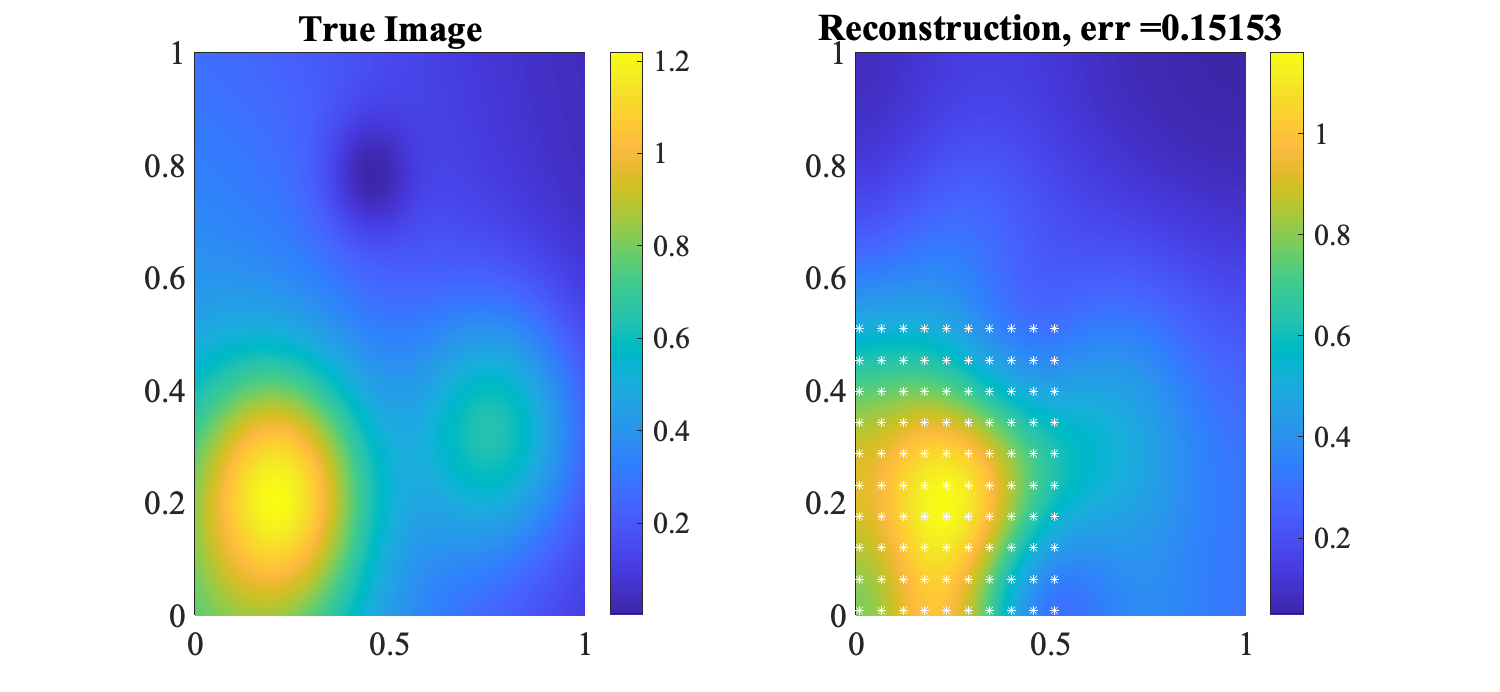} \\
    \includegraphics[scale=0.37,trim={2cm 0 2cm 0},clip]{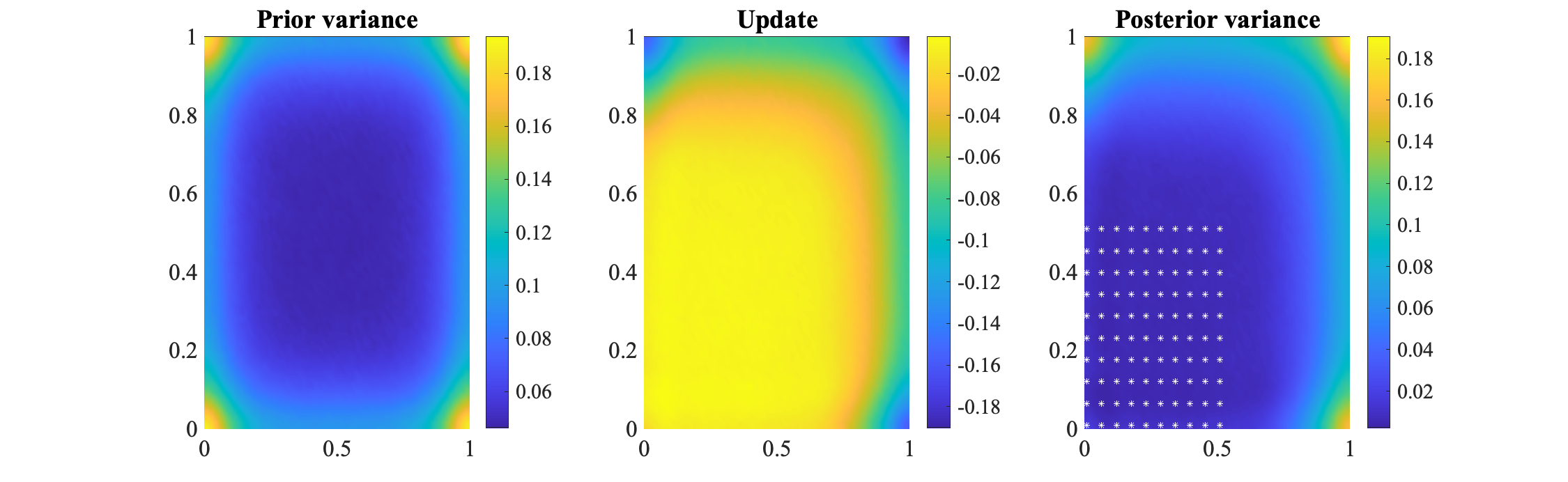}
    \caption{2D Inverse Diffusion. Top row: (left) true image, (right)  reconstruction. Bottom row: (left) prior variance $\lambda^{-2}_{\mc{C}}\bfQ$, (center) diagonals of $-\bfZ_k\bfDelta_k\bfZ_k^T$, (right) approximate posterior variance. }
    \label{fig:diffrecon}
\end{figure}

The reconstructions were performed using GenHyBR as described in Section~\ref{ssec:genhybr} and the posterior variance was approximated using the technique in Section~\ref{ssec:postvar}.  The diagonals of $\bfQ$ were estimated using the Diag++ method~\cite[Algorithm 1]{baston2022stochastic} with $300$ samples. {GenHyBR} converged in $42$ iterations and it produced a relative error of $0.15$. However, we used a basis size of $k=75$ for the subsequent uncertainty estimates.    Figure~\ref{fig:diffrecon} displays the true field and the reconstruction (top row).  The white marks indicate the sensor locations at which data is collected. In the bottom row of the same figure, we plot the prior variance (diagonals of $\lambda^{-2}_{\mc{C}}\bfQ = \bfC_{\rm pr}\bfM^{-1}$), the update $\bfZ_k\bfDelta_k\bfZ_k^T$, and the approximate posterior variance, i.e., diagonals of $\widehat{\bfGamma}_{\rm post}\bfM^{-1} = \lambda^{-2}_{\mc{C}}\bfQ - \bfZ_k\bfDelta_k\bfZ_k^T$. We see that the prior variance is high at the four corners that shows the effect of the boundary conditions. The (approximate) posterior variance shows the reduction in the uncertainty, which is especially prominent in and around the sensor coverage. 

\section{Conclusions and discussion}
In this paper, we presented efficient methods for Whittle--Mat\'ern Gaussian priors in the context of Bayesian inverse problems. However, the techniques developed here may be of interest beyond inverse problems in fractional PDEs and Gaussian random fields. Several extensions of our work are possible. First, while we focused the derivation and numerical experiments on zero Neumann boundary conditions, this is not a limitation of the framework, and it is easy to extend the framework to zero Dirichlet and Robin boundary conditions. The latter may be particularly suitable to mitigate the effects of boundary conditions (see~\cite{daon2018mitigating}). Second, the techniques in this paper may be extended to generalized Mat\'ern fields on compact Riemannian fields, see~\cite{lindgren2011explicit,lindgren2022spde,simpson2008krylov,lang2021galerkin}. Finally, in Bayesian inverse problems, it is straightforward to extend the techniques to nonlinear forward problems within a Newton-based solver for the MAP estimate~\cite{villa2018hippylib}.  Another interesting line of research is to extend these techniques to dynamic inverse problems, in which the parameters of interest change in time and we need to use spatiotemporal priors. 
\section{Acknowledgements} A.K.S. would like to thank Daniel Szyld, Daniel Sanz-Alonso, Georg Stadler, and Alen Alexanderian for helpful discussions.

\bibliography{refs}

\begin{thebibliography}{10}

\bibitem{HAntil_SBartels_2017a}
H.~Antil and S.~Bartels.
\newblock Spectral {A}pproximation of {F}ractional {PDE}s in {I}mage
  {P}rocessing and {P}hase {F}ield {M}odeling.
\newblock {\em Comput. Methods Appl. Math.}, 17(4):661--678, 2017.

\bibitem{antil2020bilevel}
H.~Antil, Z.~W. Di, and R.~Khatri.
\newblock Bilevel optimization, deep learning and fractional {L}aplacian
  regularization with applications in tomography.
\newblock {\em Inverse Probl.}, 36(6):064001, 2020.

\bibitem{HAntil_JPfefferer_SRogovs_2018a}
H.~Antil, J.~Pfefferer, and S.~Rogovs.
\newblock Fractional operators with inhomogeneous boundary conditions:
  analysis, control, and discretization.
\newblock {\em Commun. Math. Sci.}, 16(5):1395--1426, 2018.

\bibitem{bakhos2017multipreconditioned}
T.~Bakhos, P.~K. Kitanidis, S.~Ladenheim, A.~K. Saibaba, and D.~B. Szyld.
\newblock Multipreconditioned {GMRES} for shifted systems.
\newblock {\em SIAM J. Sci. Comput.}, 39(5):S222--S247, 2017.

\bibitem{baston2022stochastic}
R.~A. Baston and Y.~Nakatsukasa.
\newblock Stochastic diagonal estimation: probabilistic bounds and an improved
  algorithm.
\newblock {\em arXiv preprint arXiv:2201.10684}, 2022.

\bibitem{bolin2020rational}
D.~Bolin and K.~Kirchner.
\newblock The rational {SPDE} approach for gaussian random fields with general
  smoothness.
\newblock {\em J Comput. Graph. Stat.}, 29(2):274--285, 2020.

\bibitem{bolin2018weak}
D.~Bolin, K.~Kirchner, and M.~Kov{\'a}cs.
\newblock Weak convergence of {G}alerkin approximations for fractional elliptic
  stochastic {PDE}s with spatial white noise.
\newblock {\em BIT Numer. Math.}, 58(4):881--906, 2018.

\bibitem{bolin2020numerical}
D.~Bolin, K.~Kirchner, and M.~Kov{\'a}cs.
\newblock Numerical solution of fractional elliptic stochastic {PDE}s with
  spatial white noise.
\newblock {\em IMA J. Numer. Anal.}, 40(2):1051--1073, 2020.

\bibitem{bonito2015numerical}
A.~Bonito and J.~Pasciak.
\newblock Numerical approximation of fractional powers of elliptic operators.
\newblock {\em Math. Comput.}, 84(295):2083--2110, 2015.

\bibitem{bubba2017tomographic}
T.~A. Bubba, M.~Juvonen, J.~Lehtonen, M.~M{\"a}rz, A.~Meaney, Z.~Purisha, and
  S.~Siltanen.
\newblock Tomographic {X}-ray data of carved cheese.
\newblock {\em arXiv preprint arXiv:1705.05732}, 2017.

\bibitem{bui2013computational}
T.~Bui-Thanh, O.~Ghattas, J.~Martin, and G.~Stadler.
\newblock A computational framework for infinite-dimensional {B}ayesian inverse
  problems part {I}: The linearized case, with application to global seismic
  inversion.
\newblock {\em SIAM J. Sci. Comput.}, 35(6):A2494--A2523, 2013.

\bibitem{chung2008weighted}
J.~Chung, J.~G. Nagy, D.~P. O’leary, et~al.
\newblock A weighted {GCV} method for {L}anczos hybrid regularization.
\newblock {\em Electron. trans. numer. anal.}, 28(149-167):2008, 2008.

\bibitem{chung2017generalized}
J.~Chung and A.~K. Saibaba.
\newblock Generalized hybrid iterative methods for large-scale {B}ayesian
  inverse problems.
\newblock {\em SIAM J. Sci. Comput.}, 39(5):S24--S46, 2017.

\bibitem{chung2018efficient}
J.~Chung, A.~K. Saibaba, M.~Brown, and E.~Westman.
\newblock Efficient generalized {G}olub--{K}ahan based methods for dynamic
  inverse problems.
\newblock {\em Inverse Probl.}, 34(2):024005, 2018.

\bibitem{daon2018mitigating}
Y.~Daon and G.~Stadler.
\newblock Mitigating the influence of the boundary on {PDE}-based covariance
  operators.
\newblock {\em Inverse Probl. Imaging}, 12(5):1083, 2018.

\bibitem{flath2011fast}
H.~P. Flath, L.~C. Wilcox, V.~Ak{\c{c}}elik, J.~Hill, B.~van Bloemen~Waanders,
  and O.~Ghattas.
\newblock Fast algorithms for {B}ayesian uncertainty quantification in
  large-scale linear inverse problems based on low-rank partial {H}essian
  approximations.
\newblock {\em SIAM J. Sci. Comput.}, 33(1):407--432, 2011.

\bibitem{gazzola2019}
S.~Gazzola, P.~C. Hansen, and J.~G. Nagy.
\newblock I{R} {T}ools: a {MATLAB} package of iterative regularization methods
  and large-scale test problems.
\newblock {\em Numer. Algorithms}, 81(3):773--811, 2019.

\bibitem{herrmann2020multilevel}
L.~Herrmann, K.~Kirchner, and C.~Schwab.
\newblock Multilevel approximation of {G}aussian random fields: fast
  simulation.
\newblock {\em Math. Models Methods Appl. Sci.}, 30(1):181--223, 2020.

\bibitem{higham2008functions}
N.~J. Higham.
\newblock {\em Functions of matrices: theory and computation}.
\newblock SIAM, 2008.

\bibitem{jansson2022surface}
E.~Jansson, M.~Kov{\'a}cs, and A.~Lang.
\newblock Surface finite element approximation of spherical
  {W}hittle--{M}at\'ern {G}aussian random fields.
\newblock {\em SIAM J. Sci. Comput.}, 44(2):A825--A842, 2022.

\bibitem{kim2021hippylib}
K.-T. Kim, U.~Villa, M.~Parno, Y.~Marzouk, O.~Ghattas, and N.~Petra.
\newblock {hIPPYlib-MUQ}: A {B}ayesian inference software framework for
  integration of data with complex predictive models under uncertainty.
\newblock {\em arXiv preprint arXiv:2112.00713}, 2021.

\bibitem{lang2021galerkin}
A.~Lang and M.~Pereira.
\newblock {G}alerkin--{C}hebyshev approximation of {G}aussian random fields on
  compact {R}iemannian manifolds.
\newblock {\em arXiv preprint arXiv:2107.02667}, 2021.

\bibitem{lindgren2022spde}
F.~Lindgren, D.~Bolin, and H.~Rue.
\newblock The {SPDE} approach for {G}aussian and non-{G}aussian fields: 10
  years and still running.
\newblock {\em Spat. Stat.}, page 100599, 2022.

\bibitem{lindgren2011explicit}
F.~Lindgren, H.~Rue, and J.~Lindstr{\"o}m.
\newblock An explicit link between {G}aussian fields and {G}aussian {M}arkov
  random fields: the stochastic partial differential equation approach.
\newblock {\em J. R. Stat. Soc. Ser. B Methodol.}, 73(4):423--498, 2011.

\bibitem{roininen2014whittle}
L.~Roininen, J.~M. Huttunen, and S.~Lasanen.
\newblock {W}hittle-{M}at{\'e}rn priors for {B}ayesian statistical inversion
  with applications in electrical impedance tomography.
\newblock {\em Inverse Probl. Imaging}, 8(2):561, 2014.

\bibitem{saibaba2020efficient}
A.~K. Saibaba, J.~Chung, and K.~Petroske.
\newblock Efficient {K}rylov subspace methods for uncertainty quantification in
  large {B}ayesian linear inverse problems.
\newblock {\em Numer. Linear Algebra Appl.}, 27(5):e2325, 2020.

\bibitem{saibaba2016randomized}
A.~K. Saibaba, J.~Lee, and P.~K. Kitanidis.
\newblock Randomized algorithms for generalized {H}ermitian eigenvalue problems
  with application to computing {K}arhunen--{L}o{\`e}ve expansion.
\newblock {\em Numer. Linear Algebra Appl.}, 23(2):314--339, 2016.

\bibitem{simpson2008krylov}
D.~P. Simpson.
\newblock {\em Krylov subspace methods for approximating functions of symmetric
  positive definite matrices with applications to applied statistics and
  anomalous diffusion}.
\newblock PhD thesis, Queensland University of Technology, 2008.

\bibitem{spantini2015optimal}
A.~Spantini, A.~Solonen, T.~Cui, J.~Martin, L.~Tenorio, and Y.~Marzouk.
\newblock Optimal low-rank approximations of {B}ayesian linear inverse
  problems.
\newblock {\em SIAM J. Sci. Comput.}, 37(6):A2451--A2487, 2015.

\bibitem{stuart2010inverse}
A.~M. Stuart.
\newblock Inverse problems: a {B}ayesian perspective.
\newblock {\em Acta Numer.}, 19:451--559, 2010.

\bibitem{villa2018hippylib}
U.~Villa, N.~Petra, and O.~Ghattas.
\newblock {hIPPYlib}: {A}n extensible software framework for large-scale
  inverse problems.
\newblock {\em J. Open Source Softw.}, 3(30):940, 2018.

\bibitem{CJWeiss_BGVBWaanders_HAntil_2020a}
C.~J. Weiss, B.~G. van Bloemen~Waanders, and H.~Antil.
\newblock Fractional operators applied to geophysical electromagnetics.
\newblock {\em Geophysical Journal International}, 220(2):1242--1259, 2020.

\end{thebibliography}
\bibliographystyle{abbrv}

\end{document}